\documentclass[11pt]{article}
\usepackage{graphicx}
\usepackage{epstopdf}
\usepackage{amsfonts,amsmath, amssymb}
\usepackage{latexsym, color}
\usepackage{mathtools}
\usepackage{xcolor}
\usepackage{tikz}
\def\ds{\displaystyle}
\def\forall{\hbox{for all}~}
\def\L{{\bf L}}

\def\bfm{{\bf m}}
\def\bfv{{\bf v}}

\def\bfb{{\bf b}}
\def\bfn{{\bf n}}

\def\bfe{{\bf e}}
\def\ve{\varepsilon}
\def\n{\noindent}

\def\A{{\cal A}}

\def\inter{\hbox{\rm int}}
\def\clos{\hbox{\rm clos}}
\def\conv{\hbox{\rm conv}}
\def\avint{-\!\!\!\!\!\!\int}

\def\supp{\hbox{supp}}

\def\R{{\mathbb R}}
\def\bbS{{\mathbb S}}

\def\implies{\Longrightarrow}

\def\vs{\vskip 2em}

\def\v{\vskip 1em}

\def\dist{\hbox{\rm dist}}

\def\bega{\begin{array}}
\def\enda{\end{array}}
\def\begi{\begin{itemize}}
\def\endi{\end{itemize}}
\def\ov{\overline}

\def\Tilde{\widetilde}
\def\Hat{\widehat}

\def\N{{\mathbb N}}

\def\bel{\begin{equation}\label}
\def\eeq{\end{equation}}
\def\sqr#1#2{\vbox{\hrule height .#2pt
\hbox{\vrule width .#2pt height #1pt \kern #1pt
\vrule width .#2pt}\hrule height .#2pt }}
\def\square{\sqr74}
\def\endproof{\hphantom{MM}\hfill\llap{$\square$}\goodbreak}
\definecolor{cadmiumgreen}{rgb}{0.0, 0.42, 0.24}

\newcommand{\tarc}{\mbox{\large$\frown$}}
\newcommand{\arc}[1]{\stackrel{\tarc}{#1}}
\DeclareMathOperator*{\argmin}{arg\,min}

\newtheorem{theorem}{Theorem}[section]
\newtheorem{corollary}{Corollary}[section]
\newtheorem{lemma}{Lemma}[section]
\newtheorem{proposition}{Proposition}[section]
\newtheorem{remark}{Remark}[section]
\newtheorem{definition}{Definition}[section]

\begin{document}
\title{\bf Optimal Solutions for a Class of Set-Valued Evolution Problems}\vs

\author{Stefano Bianchini$^{(1)}$, Alberto Bressan$^{(2)}$,  and Maria Teresa Chiri$^{(3)}$
 \\~~\\
{\small $^{(1)}$~S.I.S.S.A., via Bonomea 265, Trieste 34136, Italy.}\\
 {\small $^{(2)}$~Department of Mathematics, Penn State University, University Park, Pa.~16802, USA.}\\
 {\small $^{(3)}$~Department of Mathematics and  Statistics, Queen's University,
Kingston, ON K7L3N6,
Canada.}
 \\~~\\
{\small E-mails:  bianchin@sissa.it, axb62@psu.edu, maria.chiri@queensu.ca}
}
\maketitle

\begin{abstract} 
The paper is concerned with a class of optimization problems for moving sets $t\mapsto\Omega(t)\subset\R^2$, motivated by the control of invasive biological populations. 
 Assuming that the initial contaminated set 
 $\Omega_0$ is convex, we prove that a strategy is optimal  if an only if at each given time  
 $t\in [0,T]$ the  control is active along the portion of the boundary $\partial \Omega(t)$ where the curvature is maximal.  In particular, this implies that $\Omega(t)$ is convex for all $t\geq 0$.
The proof relies  on the analysis of a one-step constrained optimization problem, obtained by a time discretization.
\end{abstract}

\section{Introduction}
\label{s:1}
\setcounter{equation}{0}

Motivated by a model in \cite{BCS1, BCS2}, describing the control of an invasive  biological 
species, we consider  here the evolution problem for a set $\Omega(t)\subset\R^2$
of finite perimeter,  depending on the normal velocity assigned at every 
boundary point $x\in \partial \Omega(t)$.
We think of $\Omega(t)$ as the contaminated set at time $t\geq 0$.  
If no control is applied, this set expands with unit speed in all directions.
By implementing a control strategy, we assume that one can reduce the area of $\Omega(t)$ at rate $M$ per unit time.

To model this situation, for $t\in [0,T]$ and $x\in \partial\Omega(t)$, 
we denote by $\beta(t,x)$ the normal speed of the boundary at the point $x$, in the direction
of the interior normal.  In other words, if the sets $\Omega(t)$ are described by
$$\Omega(t)~=~\Big\{x=(x_1,x_2)\in \R^2\,;~~\psi(t,x_1, x_2)>0\Big\}$$
for some differentiable function $\psi$,
then 
$$\beta~\doteq~{-\psi_t\over\sqrt{\psi_{x_1}^2 + \psi_{x_2}^2}}\,.$$
We denote by
\bel{Edef} E(\beta)~\doteq~\max~\{1+\beta, \, 0\}\eeq
the {\bf control effort}, needed to push the boundary of $\Omega$ inward
with speed~$\beta$.
\begin{definition}\label{d:11} Given a constant $M>0$, we say
that a set-valued function $t\mapsto \Omega(t)$ is {\bf admissible} if the corresponding  characteristic function $t\mapsto {\bf 1}_{\Omega(t)}$ is Lipschitz continuous from $[0,T]$
into $\L^1(\R^2)$, and 
moreover
\bel{eb}
\int_{\partial\Omega(t)}E\bigl(\beta(t,x)\bigr)\, \mathcal H^1(dx)~\leq~ M,\qquad \hbox{for a.e.}~~ t\in [0,T].\eeq
Here  $\beta$ denotes the velocity of a boundary point in the inward normal direction, 
and the integral is taken w.r.t.~the 1-dimensional Hausdorff measure
along the boundary of $\Omega(t)$.
\end{definition}


Given an initial set $\Omega_0\subset\R^2$ and a constant $M>0$, 
three problems will be considered.
\begi
\item[{\bf (NCP)}] {\bf Null Controllability Problem.}  
{\it  Find an admissible  set-valued function $t\mapsto \Omega(t)$  and a time $T>0$ such that}
\bel{nco}\Omega(0)~=~\Omega_0,\qquad\qquad 
\Omega(T)~=~\emptyset.\eeq
\item[{\bf (MTP)}]  {\bf Minimum Time Problem.}  {\it  Among all admissible strategies that satisfy
(\ref{nco}), find one that minimizes the time $T$.}

\item[{\bf (OP)}] {\bf Optimization Problem.}  {\it  Given a time interval $[0,T]$ and constants $c_1,c_2\geq 0$, 
find an admissible set-valued function $t\mapsto \Omega(t)$  which minimizes the cost
\bel{cost1}
J~=~
c_1\int_0^T {\mathcal L}^2 \bigl(\Omega(t)\bigr)\, dt + c_2\,{\mathcal L}^2 \bigl(\Omega(T)\bigr),\eeq
subject to $\Omega(0)=\Omega_0\,$.
}
\endi
Here and in the sequel, we use the notation
${\mathcal L}^2 (\Omega)$ to denote the 2-dimensional Lebesgue measure
of a set $\Omega\subset\R^2$, while ${\mathcal H}^1 (\partial\Omega)$ will be used 
to denote the 
1-dimensional Hausdorff measure of its boundary.

%

\begin{remark}{\rm When the control effort is zero, the set $\Omega(t)$ 
expands in all directions
with unit speed.
On the other hand,  the bound (\ref{eb}) on the instantaneous control effort allows us reduce the area of $\Omega(t)$  at rate  $M$ per unit time.  This yields a basic relation between 
the growth rate of the area of $\Omega(t)$ and its perimeter:
\bel{AS}{d\over dt} {\mathcal L}^2 \bigl(\Omega(t)\bigr)~=~{\mathcal H}^1 \bigl(\partial\Omega(t)\bigr) -M.\eeq
The Null Controllability Problem {\bf (NCP)} can thus be solved if and only if we can reduce the perimeter $P(t)={\mathcal H}^1 \bigl(\partial\Omega(t)\bigr)$
to a value strictly smaller than $M$.
}
\end{remark}

A more general class of optimization problems for moving sets was recently considered in 
\cite{BCS2}, proving the existence of optimal strategies and deriving some necessary conditions for optimality.

In this paper  we consider the 
optimization problems {\bf (MTP)} and {\bf (OP)}, assuming that the initial set $\Omega_0$ is convex.
Our main result, Theorem~\ref{t:61}, completely characterizes the optimal strategies.  Confirming a conjecture
proposed in \cite{BCS2}, we prove that
a strategy is optimal if an only if, at each given 
time  $t\in [0,T]$, the 
control is active precisely along the portion of the boundary $\partial \Omega(t)$ 
where the curvature is maximal.  We observe that, with this control,  the 
perimeter $P(t)$ also shrinks at the fastest possible rate.  Moreover, all sets $\Omega(t)$ remain convex.

As a preliminary to the proof of the main theorem, Section~\ref{s:2} collects 
several geometric results concerning $r$-semiconvex sets, i.e., sets that 
satisfy the outer sphere condition with radius $r$.  In particular, we prove a sharp bound on their perimeter, and give estimates on how the area of their $r$-neighborhood
changes, when the boundary is perturbed.

In Section~\ref{s:3} we study a one-step minimization problem, derived from 
the original evolution problem by discretizing time. More precisely, given a compact convex set $U\subset\R^2$ and a constant $0<a<{\cal L}^2(U)$, we seek a subset
$\Omega\subset U$ with area ${\cal L}^2(\Omega)=a$, such that the area of 
its $r$-neighborhood  $B_r(\Omega)$ is as small as possible.
By a  detailed analysis, we prove that this optimal set $\Omega$ is always convex, see 
Theorem~\ref{t:31}.   
As a consequence, the set  $\Omega$ is also optimal for the problem 
of minimizing the perimeter, subject to the same constraints. 
In the literature, this second problem 
has already been 
studied in \cite{SZ}, and is well understood  in dimension $n=2$.
 Properties of constrained perimeter-minimizing
sets, described in Section~\ref{s:4},  play a key role in determining the optimal strategy $t\mapsto\Omega(t)$ 
for both {\bf (OP)} and {\bf (MTP)}.


For an introduction to geometric measure theory and BV functions we refer to \cite{AFP, EG, M}.
Various other models of moving sets, subject to external control, have been considered in 
\cite{Bblock, Breview, BMN, BZ, CPo, CLP}.  More detailed models of the control of an invasive
biological species can be found in \cite{ACD, ACM}. 
See also \cite{RBZ} for 
related results on controlled reaction-diffusion equations.

\section{Preliminary geometric lemmas}
\label{s:2}
\setcounter{equation}{0}

Throughout the following,
$B_r(x)$ denotes the open ball centered at $x$ with radius $r$, while
\begin{equation*}
\Omega^r ~\doteq~ \Omega + B_r(0) ~=~ \big\{x \,;~ \hbox{\rm dist}(x,\Omega) < r \big\}
\end{equation*}
denotes the open neighborhood of radius $r$ around the set $\Omega\subset\R^2$,
{and we will use the standard notation $\mathcal L^2$ for the Lebesgue measure in $\R^2$ and $\mathcal H^1$ for the $1$-dimensional Hausdorff measure.} 
The closure, the interior, and the convex hull of a set $\Omega$ are denoted by 
$\clos\, \Omega$, $\inter\,\Omega$ and $\conv\, \Omega$, respectively.
Given a vector $\bfv=(v_1, v_2)\in\R^2$, we write $\bfv^\perp = (-v_2, v_1)$ for the perpendicular vector.

\begin{figure}[ht]
\centerline{\hbox{\includegraphics[width=13cm]{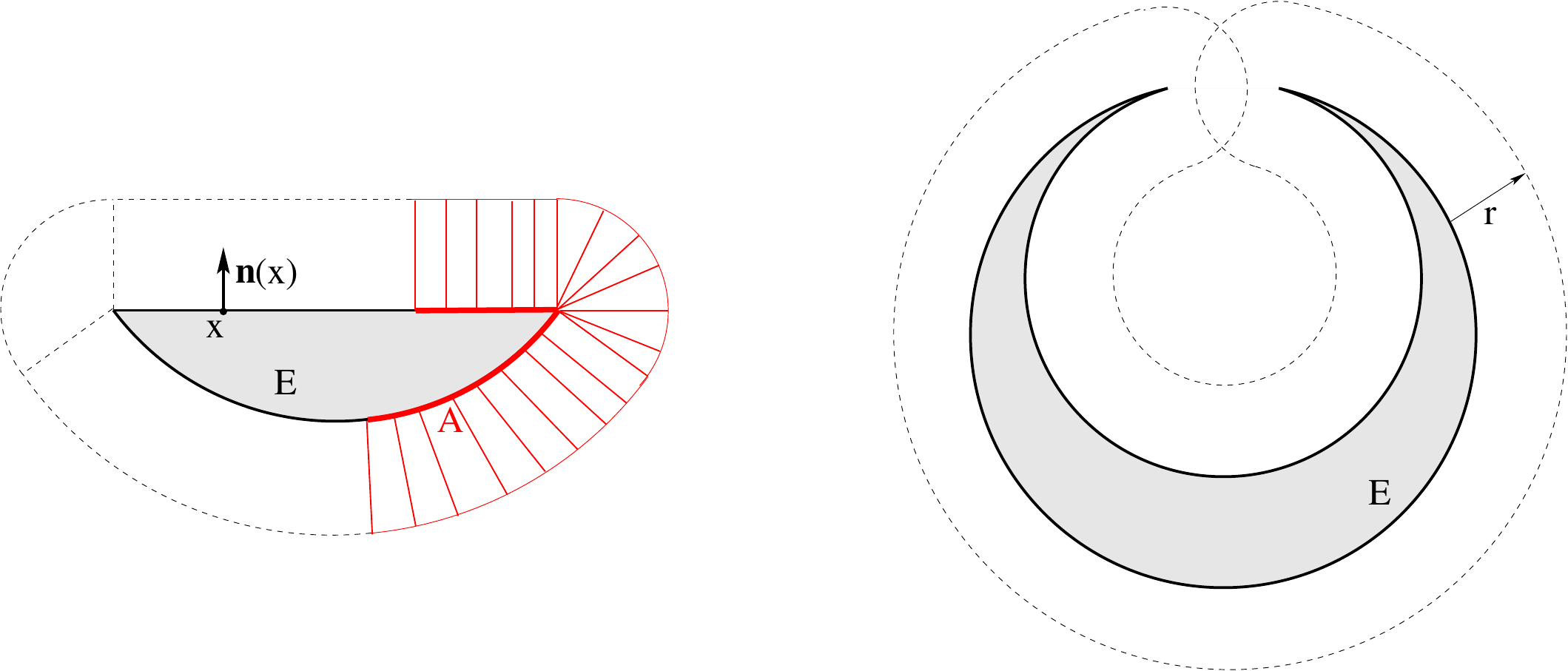}}}
\caption{\small   Left: a convex set $E$ and its  $r$-neighborhood. Given a measurable subset of the
boundary $A\subset \partial E$,  the set of points $y\in E^r$ that project onto $A$
has area computed by (\ref{33}).    Right: if the set $E$ is $r$-convex but not convex, the
formula (\ref{36}) may hold only as an inequality, because of the overlap.}
\label{f:sc85}
\end{figure}

Let $E\subset\R^2$ be a compact convex set with nonempty interior, and let  $t\mapsto \gamma(t)$ be an arc-length parametrization of the boundary $\partial E$, oriented counterclockwise.
Call $\bbS^1= \{\bfe\in \R^2\,;~|\bfe|=1\}$ the set of unit vectors in $\R^2$.
For every boundary point $x\in \partial E$ consider the set of outer unit normals
\begin{equation*}
\bfn(x) ~\doteq~ \bigl\{\mathbf e \in \bbS^1\,;~~\bfe \cdot (x-y) \geq 0 \quad\forall y \in E \bigr\}.
\end{equation*}
We observe that 
the (possibly multivalued) map   $t\mapsto \bfn(t)\doteq \bfn\bigl(\gamma(t)\bigr)\subset \bbS^1$
is a set-valued function with closed graph and connected values.
We will sometimes use the notation $\bfn(t) = e^{i\theta(t)}$. In this case $t \mapsto \theta(t)$ is a monotone multifunction. 
A few observations are in order.
\v
{\bf 1.} The map $t\mapsto \bfn(t)$ is a BV function, with total variation $TV(\bfn(\cdot )) = 2\pi$.
The absolutely continuous and the singular part of its measure-valued derivative
will be denoted by
$$
D\bfn~ = ~D^{ac} \bfn + D^{sing} \bfn ~ \in~  \mathcal M(\partial \Omega,\bbS^1).
$$
Note that with this notation we require that in the jump points the quantity $\bigl|\bfn(x+) - \bfn(x-)\bigr|$ is the length of the smaller arc $\bigl[\bfn(x-),\bfn(x+)\bigr]$, i.e. the jump $\bigl|\theta(t+) - \theta(t-)\bigr|$, 
$x = \gamma(t)$.
\v
{\bf 2.}  If $A \subset \partial E$ is a Borel subset, then (see Fig.~\ref{f:sc85}, left)
\bel{32}
\mathcal H^1 \bigg(  \bigcup_{x \in A} x + h \bfn(x) \bigg) ~=~ \mathcal H^1(A) + h |D\bfn|(A) ~= ~\mathcal H^1(A) + h D\theta(A),
\eeq
\bel{33}
\mathcal L^2 \bigg( \bigcup_{x \in A, ~\rho\in [0,h]} x + \rho \bfn(x) \bigg) ~=~ h \mathcal H^1(A) + \frac{1}{2} h^2 |D\bfn|(A)~ =~ h \mathcal H^1(A) + \frac{1}{2} h^2 D\theta(A).\eeq
These formulas generalize the classical Steiner's formulas for the
perimeter and the area of the $r$-neighborhood around a convex set
(see Theorem 10.1 in \cite{G}).   They can be proved by approximating $E$ with a polygon and passing to the limit.
Notice that the second is the integral of the first one, by the coarea formula.
The above identities can be extended to a more general class of sets with finite perimeter.

\begin{definition}\label{d:dual}
Given a closed set $E\subset \R^2$, an open set $F$ and a radius $r>0$, we  write
\begin{equation*}
E^r \,\doteq\,  \bigcup_{x \in E} B_r(x), \qquad\qquad
 F^{-r}\,\doteq\,    \bigg( F \setminus \bigcup_{y \notin F} B_r(y) \bigg) ~=~ \R^2 \setminus (\R^2 \setminus F)^r . 
\end{equation*}
We say that $E,F$ are in $r$-duality (or simply in duality) if $F = E^r$ and  $E = F^{-r}$.
\end{definition}

From the above definition, it
immediately follows that the two sets $E^r$ and $(E^r)^{-r}$ are in duality.
Note that we always have $E\subseteq (E^r)^{-r}$, but equality does not hold, in general.

\begin{definition} Let $E\subset\R^2$ be a closed set, and let $r>0$.
We say that $E$ has the {\bf exterior $r$-ball property}, or equivalently that $E$ is {\bf $r$-semiconvex}, if 
for every boundary point $x \in \partial E$ there is an outer ball $B_r(y) \subset \R^2 \setminus E$ with 
$|x-y| = r$.  When this holds, we say that the segment with endpoints $x,y$ is an {\bf optimal ray}.
\\
We say that an open set $F$ has the  {\bf interior $r$-ball property}, or equivalently that $F$ 
is {\bf $r$-semiconcave}, if
every point $x\in F$ is contained in some open ball $B_r(y)\subseteq F$.  
\end{definition}

In the following, having fixed the radius $r$, for shortness we will just say semiconvex/semiconcave.
As immediate consequences of  the above definitions, one has:
\begi\item[(i)] $E$ is $r$-semiconvex iff $E = (E^r)^{-r}$.
\item[(ii)] $F$ is $r$-semiconcave iff $F= (F^{-r})^r$.
\item[(iii)] For every closed set $E$, 
every point $x\in \partial E^r$ belongs to the boundary of a ball of radius $r$ contained in $E^r$, with center in $E$;
\item[(iv)]  For every closed set $E$, every point  $x\in \partial (E^r)^{-r}$ belongs to the boundary of a ball of radius $r$ contained in $\R^2 \setminus E$, with center in $\partial E^r$.
\endi

We observe that, if $E$ is compact, $r$-semiconvex with $\inter\, E$ connected but not simply connected, then each of its ``holes" must contain an open ball of radius $r$.  Therefore, the complement $\R^2\setminus E$ can have at most finitely many 
connected components.  To fix ideas, we call $V_1$ the unbounded component, 
and $V_2,\ldots, V_N$ the bounded components.  Each boundary $\partial V_k$, $k=1,\ldots,N$
is a simple closed curve with finite length.   Let $t\mapsto \gamma_k(t)$ be an arc-length
parameterization of $\partial V_k$, oriented counterclockwise in the case of $V_1$ and 
clockwise for $V_2,\ldots, V_N$.
As before, let $\bfn(x)$ be the set of outer unit normals at the point $x\in \partial E$.
Using complex notation, we again write $\bfn_k(t)=\bfn(\gamma_k(t))=e^{i\theta_k(t)}$
for the set of  unit outer normal vectors at the point $\gamma_k(t)$.  The assumption of $r$-semiconvexity
implies that the (possibly multivalued) map $t\mapsto \theta_k(t)$ has bounded variation.
Its distributional derivative satisfies
\bel{Dthe}
D\theta_k ~\geq~ -\frac{1}{r} \mathcal L^1.
\eeq
Indeed, the negative part of the measure $D\theta$ is absolutely continuous and has
uniformly bounded density w.r.t.~1-dimensional Lebesgue measure ${\mathcal L^1}$ on the unit
circumference $\bbS^1$.  For future use, we 
notice that this yields:
\begin{lemma}\label{l:31} 

\begi \item[]
\item[(i)] The set 
$\R^2\setminus V_1$ is not convex if and only if  there exists a  point $x^*= \gamma_1(t^*)$
where the function $t\mapsto \theta_1(t)$ is differentiable, with a strictly negative derivative.
\item[(ii)] For every $k=2,\ldots,N$, there exists  a point $x^*= \gamma_k(t^*)$
where the function $t\mapsto \theta_k(t)$ is differentiable, with a strictly negative derivative.
\endi
\end{lemma}

The identities (\ref{32})-(\ref{33}) have counterparts for general $r$-semiconvex sets. More precisely,
consider a Borel subset $A\subseteq \partial V_k\subseteq \partial E$. For $0 < h \leq r$ we then have
\bel{35}
\mathcal H^1 \bigg( \bigcup_{x \in A} x + h \bfn(x) \bigg) ~\leq~ \mathcal H^1(A) + h D\theta_k(A),
\eeq
\bel{36}
\mathcal L^2 \bigg(  \bigcup_{x \in A, ~\rho\in [0,h]} x + \rho \bfn(x) \bigg) ~\leq~ h \mathcal H^1(A) + \frac{1}{2} h^2 D\theta_k(A).
\eeq
These formulas can again be  proved by approximating the set
$E$ with polygons, then passing to the limit.  The second one is obtained 
from the first by an integration.   We observe that, if $E$ is not convex, then 
there can be distinct points $x,y\in \partial E$ such that 
the sets $\bigl\{x+ \rho\bfn(x)\,;~\rho\in [0, h]\bigr\} $, $\bigl\{y+ \rho\bfn(y)\,;~\rho\in [0, h]\bigr\} $
have non-empty intersection (see Fig.~\ref{f:sc85}, right).  This motivates the inequality signs in (\ref{35})-(\ref{36}).

\subsection{The perimeter of a semiconvex set.}

\begin{lemma} \label{l:22} Any $r$-semiconvex set $E$ has locally finite perimeter. Indeed, for any ball of radius $r$,
one has
\bel{brp}
{\mathcal H}^1 \bigl( B_r(\bar x) \cap \partial E\bigr)~ \leq~~2 \pi r.\eeq
\end{lemma}

{\bf Proof.} {\bf 1.} As a first  step, we observe that 
the boundary $\partial E$  is rectifiable with length locally finite. Indeed,  for any unit vector 
 $\bfe\in\bbS^1$, consider the set of points  $x\in \partial E$ which have an optimal ray of direction $\bfn(x)$ such that $\bigl|\bfn(x)-\bfe\bigr|<1/2$. By the outer ball property of $r$-semiconvex sets, for every such an $x$ there is an open cone $C_x$ with axis in the direction $\pm \bfe$ and opening $>\pi/6$ such that $\partial E \cap C_x \cap B_r(x) = \emptyset$: indeed, $C_x \cap B_r(x) \subset B_r(x + r \bfe) \subset \R^2 \setminus E$.
A standard rectifiability criterion (see Theorem 2.61 in \cite{AFP})  shows that $\partial E$ is rectifiable, and moreover the fact that the curves are separated gives that their total length is locally finite.
\v
{\bf 2.} 
Next, we claim that, for any fixed ball $B_r(\bar x)$ and any  finite family of balls $B_r(x_i)$, $i=1,\dots,n$, there holds
\begin{equation}
\label{ede}
{\mathcal H}^1 \bigg( B_r(\bar x) \cap \partial \bigcup_{i=1}^n B_r(x_i) \bigg)~ \leq~
 {\mathcal H^1} \bigg( \partial B_r(\bar x) \cap \bigcup_{i=1}^n B_r(x_i) \bigg)~ \leq ~2 \pi r.
\end{equation}
Indeed,
the above estimate is trivial when $n=1$.  
By induction, assume that it holds for $n$ and consider an additional  ball $B_r(x_{n+1})$.
By relabeling,  we can assume that $\partial B_r(x_{n+1}) \setminus \cup_i^n \clos\,B_r(x_i)$ 
is an arc intersecting $\partial B_r(\bar x)$ in at least one point. 

\begin{figure}[ht]
\centerline{\hbox{\includegraphics[width=14cm]{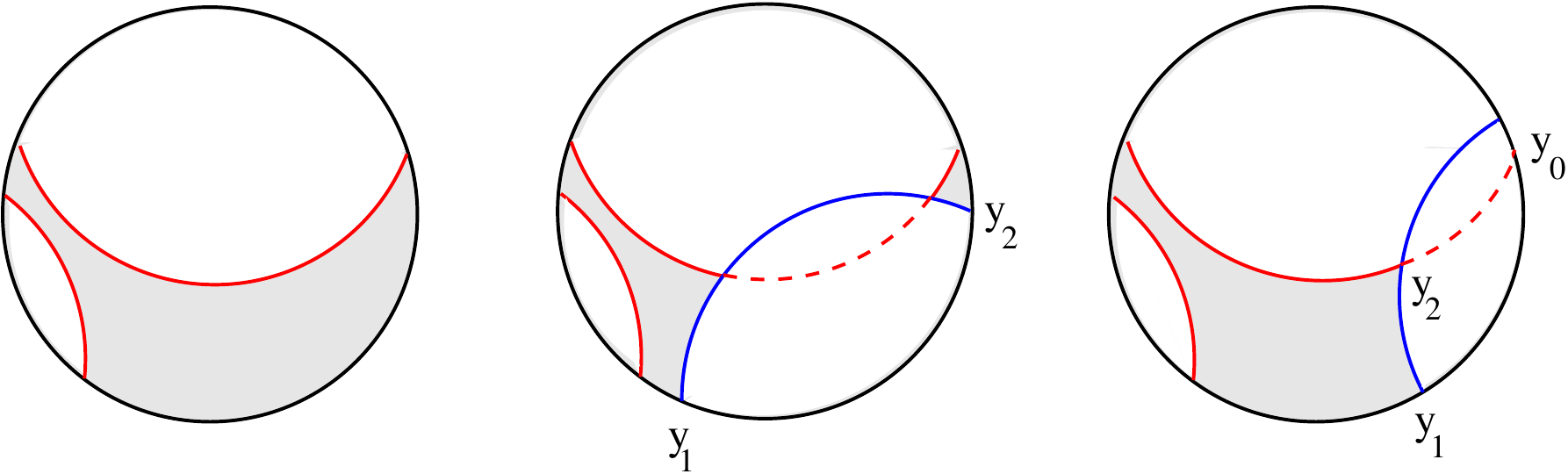}}}
\caption{\small  Left: the set $B_r(\bar x) \setminus \bigcup_{i=1}^n B_r(x_i)$.
Center and right:  the set $B_r(\bar x) \setminus \bigcup_{i=1}^{n+1} B_r(x_i)$, in the two cases
considered in step {\bf 2} of the proof of Lemma~\ref{l:22}.}
\label{f:sc106}
\end{figure}

\begi
\item
If the intersection 
\bel{ints}\partial B_r(x_{n+1}) \cap \left( \partial B_r(\bar x)  \setminus \bigcup_{i=1}^n B_r(x_i)\right)\eeq
consists of two distinct points $y_1,y_2$
(see Fig.~\ref{f:sc106}, center), then either \\
(i) there is a ball $B_r(x_i)$ such that $B_r(x_i) \cap B_r(\bar x) \subset B_r(x_{n+1}) \cap B_r(\bar x)$. In this case the ball $B_r(x_i)$ can be removed from the family of balls and we can apply the recurrence hypothesis.
Or  else \\
(ii)  we observe that 
\bel{2pt2}\bega{lr}{\mathcal H}^1 \bigg( B_r(\bar x) \cap \partial \bigcup_{i=1}^{n+1} B_r(x_i) \bigg)~\\=~{\mathcal H}^1 \bigg( B_r(\bar x) \cap\Big[\Big(\partial B_r(x_{n+1})\setminus\bigcup_{i=1}^{n} B_r(x_i)\Big)\cup \Big(\partial\bigcup_{i=1}^{n} B_r(x_i)\setminus B_r(x_{n+1}) \Big)\Big]\bigg) \\
\leq~ {\mathcal H}^1 \bigg( B_r(\bar x) \cap \partial B_r(x_{n+1})\bigg)\,+\, {\mathcal H}^1 \bigg( B_r(\bar x) \cap \partial\bigcup_{i=1}^{n} B_r(x_i)\bigg)\\
\leq~ {\mathcal H}^1 \bigg(  \partial B_r(\bar x) \cap B_r(x_{n+1})\bigg)\,+\, {\mathcal H}^1 \bigg(\partial B_r(\bar x) \cap\bigcup_{i=1}^{n} B_r(x_i)\bigg)\\
=~ {\mathcal H}^1 \bigg(\partial B_r(\bar x) \cap  \bigcup_{i=1}^{n+1} B_r(x_i) \bigg).\enda\eeq
Hence the first inequality in 
(\ref{ede})  holds. 
  \item 
On the other hand,  if the intersection (\ref{ints}) consists of a single point $y_1$ (see Fig.~\ref{f:sc106}, right), 
then   it contributes to the measure above by an arc 
${ \arc{y_1y_2}} \subset \partial B_r(x_{n+1}) \cap B_r(\bar x)$. In this case, the arc ${\arc{y_1y_2}}$ 
replaces the set $B_r(x_{n+1}) \cap B_r(\bar x) \cap \partial \Big(\bigcup_{i=1}^n B_r(x_i)\Big)$, 
and either\\
(i)  there is ball $B_r(x_i)$ which is not contributing to 
$B_r(\bar x) \cap \partial \bigl(\bigcup_{i=1}^{n+1}\textflorin B_r(x_i)\bigr)$. 
This happens when two or more balls are involved in the set 
$B_r(x_{n+1}) \cap B_r(\bar x) \cap \partial \Big(\bigcup_{i=1}^n B_r(x_i)\Big)$.  
In this case, the ball $B_r(x_i)$ can be removed, and the result follows by the inductive assumption. Or else 
\\
(ii) the arc ${\arc{y_1y_2}}$ 
replaces an arc ${\arc{y_0y_2}} = B_r(x_{n+1}) \cap \partial B_r(x_j)$. 
A simple computation shows that
\begin{equation*}
\text{length of ${\arc{y_1y_2}}$}~ \leq ~\text{length of ${\arc{y_0y_2}}$} + \text{length of ${\arc{y_0y_1}}$},
\end{equation*}
where ${\arc{y_0y_1}}$ is the the arc of $\partial B_r(\bar x)$ which belongs to $B_r(x_{n+1}) \cap B_r(\bar x) \setminus \cup_i^n B_r(x_i)$. Hence again \eqref{ede}  is verified.
\endi
By induction, the estimate (\ref{ede}) holds for every $n\geq 1$.

{\bf 3.}  To  complete the proof, we need to consider 
the general case where $E$ is any $r$-semiconvex set. 
Since $\partial E \cap B_r(\bar x)$ is rectifiable with finite length, by the 
Vitali-Besicovitch Covering Theorem { \cite[Theorem 2.19]{AFP}} for every $\ve > 0$ there is a finite set of points $x_i \in \partial E \cap B_r(\bar x)$ and radii $r_i < \ve$, $i=1,\dots,N$, such that
\begin{enumerate}
\item $x_i$ belongs to the reduced boundary $\partial^* E$, i.e. it has a unique normal $\bfn_i$;
\item $B_{r_i}(x_i) \subset B_r(\bar x)$;
\item there is a Lipschitz curve $\gamma_i$ such that
\[
|\dot \gamma_i - \bfn_i^\perp| < \ve \quad \text{and} \quad \mathcal H^1 \big( \mathrm{graph}(\gamma_i) \cap B_{r_i}(x_i) \big) > (1 - \ve) 2r_i; 
\]
\item it holds
\begin{equation*}
\sum_i \big| \mathcal H^1(\partial E \cap B_r(x_i)) - 2 r_i \big| < \ve, \qquad \bigg| \mathcal H^1(\partial E \cap B_r(\bar x)) - \sum_i 2r_i \bigg| < \ve.
\end{equation*}

\end{enumerate}
Indeed, the balls satisfying the first 3 statements are a fine 
covering of $\partial^* E \cap B_r(\bar x)$.

\begin{figure}
\centering{\resizebox{8cm}{!}{\includegraphics{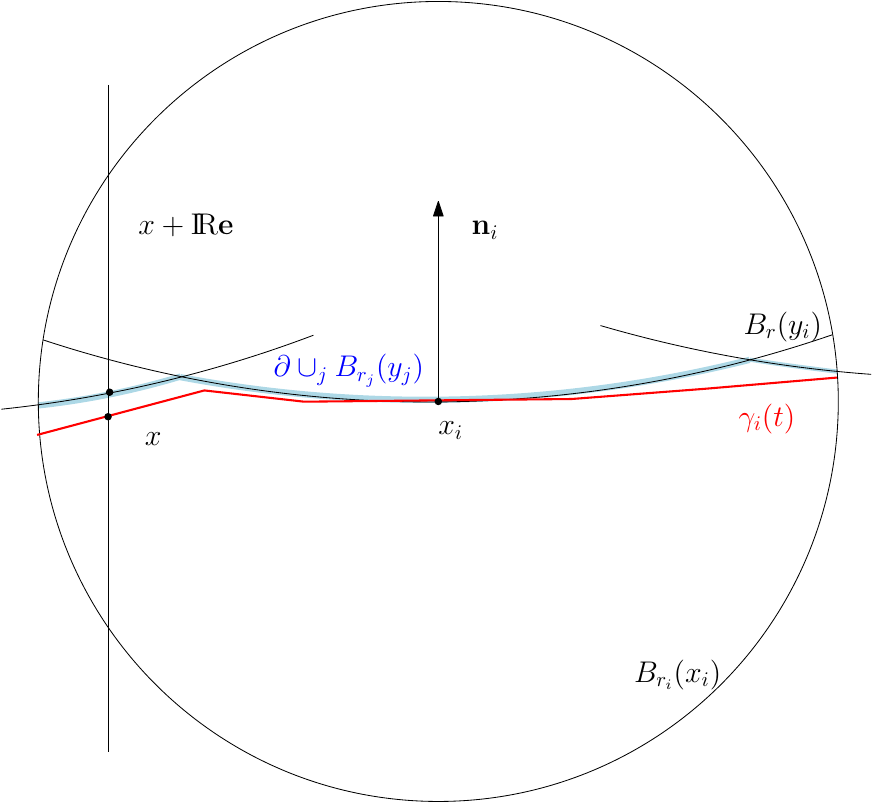}}}
\caption{Up to two small end parts, to each point of $x \in \partial E \cap \mathrm{graph}(\gamma_i) \cap B_{r_i}(x_i)$ (red) there corresponds at least one point of $\partial \cup_j B_r(y_j)$ (light blue), giving (\ref{Equa:approx_boundary_est}).}
\label{Fig:approx_boundary}
\end{figure}

For each $x_i$, consider the optimal ray $[x_i,y_i]$ (there is only one because $x_i \in \partial^* E$) and the family of balls $B_r(y_i) \subset \R^2 \setminus E$. The boundary of the set $\cup_i B_r(y_i)$ 
consists of finitely many curves, and inside every ball $B_{r_i}(x_i)$ its length must be at least
\begin{equation}
\label{Equa:approx_boundary_est}
\mathcal H^1 \bigg( B_{r_i}(x_i) \cap \partial \bigcup_j B_r(y_j) \bigg) ~>~ (1 - 2\ve) (2 r_i)
\end{equation}
for $r_i \ll 1$. Indeed, for each point $x \in \partial E \cap \mathrm{graph}(\gamma_i) \cap B_{(1-\ve)r_i}(x_i)$ the line $x + \R \bfn_i$ must intersect the boundary of $\cup_i B_r(x_i)$ at one point inside $B_{r_i}(x_i)$ (see Fig. \ref{Fig:approx_boundary}).

 In view of (\ref{ede}),  we thus conclude that
\begin{equation*}
\begin{split}
\mathcal H^1(\partial E \cap B_r(\bar x)) &\leq ~\ve +  \sum_i 2r_i ~\leq ~
\ve + \frac{1}{1-2\ve} \sum_i \mathcal H^1 \bigg( B_{r_i}(x_i) \cap \partial \bigcup_j B_r(y_j) \bigg) \\
&\leq ~\ve + \frac{1}{1 - 2 \ve} \mathcal H^1 \bigg( \partial B_r(\bar x) \cap \bigcup_j B_r(y_j) \bigg)
~ < ~\ve + \frac{2\pi r}{1 - 2 \ve}.
\end{split}
\end{equation*}

Taking the limit $\ve \searrow 0$ we obtain the bound (\ref{brp}).
%
\endproof

\begin{remark} {\rm
In the case of different radii  $0<\rho<r$,  the same arguments used in the 
proof of Lemma~\ref{l:22}
show that the
estimate (\ref{ede}) can be replaced by 
\bel{edf}
{\mathcal H^1} \bigg( B_\rho(\bar x) \cap \partial \bigcup_{i=1}^n B_r(x_i) \bigg)~ \leq ~ 2 \pi \rho. 
\eeq }
\end{remark}

\subsection{A property of $r$-semiconvex sets.}
Throughout this section we consider a compact set $E\subset \R^2$ whose boundary is 
a simple closed curve $s\mapsto \gamma(s)$, parameterized by arc-length and oriented counterclockwise. As before, the unit outer normals are denoted by
$e^{i\theta(s)} = \bfn(s) = - \dot \gamma(s)^\perp$. 

If the set $E$ is $r$-semiconvex, then the negative part of of the distributional derivative of 
$\theta$ is absolutely continuous w.r.t.~one-dimensional Lebesgue measure. Namely,  
 \bel{tsc}D \theta~ \geq ~- {1\over r} {\mathcal L}^1.\eeq
 
\begin{figure}[ht]
\centerline{\hbox{\includegraphics[width=5cm]{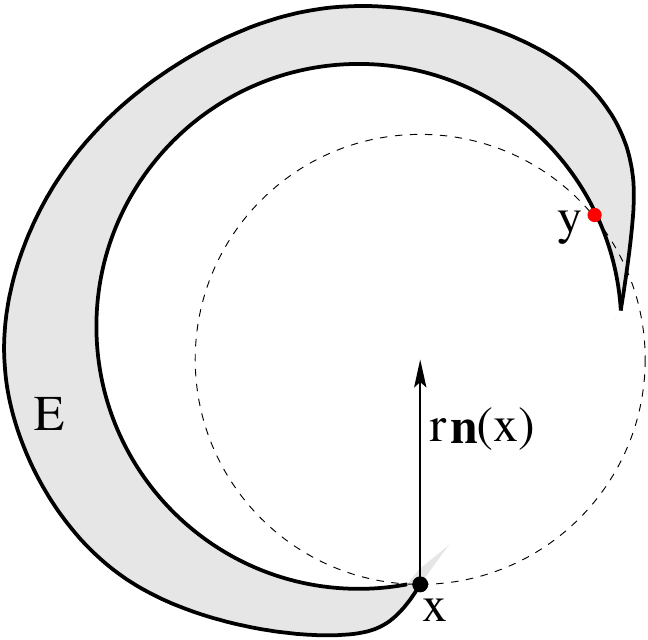}}\qquad\qquad\qquad\hbox{\includegraphics[width=6cm]{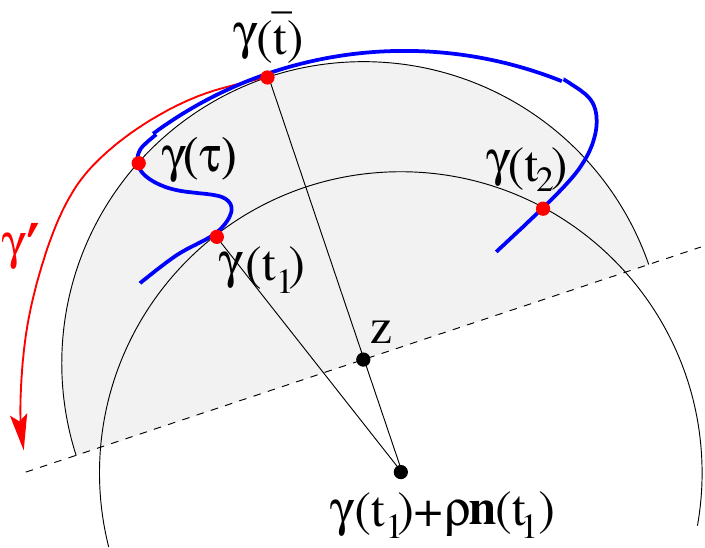}}}
\caption{\small Left: if (\ref{tsc}) holds, the outer curvature  radius is $\geq r$.  Hence in a 
neighborhood of $x$ the set $E$ lies outside the ball $B_r\bigl(x+ r\bfn(x)\bigr)$. However, the boundary $\partial E$ can enter this ball at a point $y$, where the boundary arc ${\arc{xy}}$ has length $>\pi r$.  Right: the points $\gamma(t_1)$, $\gamma(t_2)$, $\gamma(\bar t)$ and 
$\gamma(\tau)$, considered in step {\bf 2}
of the proof of Lemma~\ref{l:23}. Starting from $\gamma(\bar t)$ and moving toward 
$\gamma(t_1)$, we enter the ball $B_\rho(z)$ (shaded region) at a point $\gamma(\tau)$,
before reaching $\gamma(t_1)$.     Indeed, any curve $\gamma'$ of length
$\leq \pi\rho/2$ starting at $\gamma(\bar t)$ and remaining outside the ball $B_\rho(z)$ 
cannot touch the ball $B_\rho\bigl(\gamma(t_1) + \rho\bfn(t_1)\bigr)$.
The inductive process eventually identifies a point 
$\gamma(t^*)$ where the outer radius of curvature is $\leq\rho<r$, thus obtaining a contradiction.
  }
\label{f:sc107}
\end{figure}

Conversely,  assume that  $e^{i\theta} = - \dot \gamma^\perp$ and that (\ref{tsc}) holds, so that
  $\theta$ has bounded variation.  One can then define the set of unit normal vectors as
\begin{equation*}
\bfn(t) = \big\{e^{i\phi}, \phi \in [\theta(t-),\theta(t+)] \big\}.
\end{equation*}
Equivalently, $\bfn(t)$ is the set of unit vectors $\bfn\in\bbS^1$  such that
\begin{equation*}
\limsup_{\delta \searrow 0} \frac{\sup \big\{\bfn \cdot (y - \gamma(t))\,;\quad  y \in E \cap B_\delta(\gamma(t)) \big\}}{\delta} ~\leq ~0.
\end{equation*}
Let (\ref{tsc}) hold and consider any boundary point $x\in\partial E$. As shown in Fig.~\ref{f:sc107}, by
moving along the boundary of $E$,  it is possible to get into 
the interior of the outer tangent ball 
$B_r\bigl(x+r\bfn(x)\bigr)$, but only after having travelled along an arc of length $>\pi r$.

\begin{lemma}\label{l:23} In the above setting, if (\ref{tsc}) holds, then  for every $x = \gamma(\bar t) \in \partial E$ one has
\bel{eter}
B_r\bigl(\gamma(\bar t) + r \bfn(\bar t) \bigr) \cap \Big\{\gamma(t)\,;\quad 
 t \in [\bar t -\pi r\,,~ \bar t +\pi r] \Big\} ~= ~\emptyset. 
\eeq
\end{lemma}
{\bf Proof.}  If (\ref{eter}) fails, then there exists $0<\rho<r$ such that 
\bel{etero}
B_r\bigl(\gamma(\bar t) + \rho \bfn(\bar t) \bigr) \cap \Big\{\gamma(t)\,;\quad 
 t \in [\bar t -\pi \rho\,,~ \bar t +\pi \rho]  \Big\} ~\not= ~\emptyset. 
\eeq
A contradiction can then be obtained in two steps.

\v
{\bf 1.} We claim that, for every $\bar t$, there exists $\ve>0$ such that 
\bel{etee}
B_r\bigl(\gamma(\bar t) + r \bfn(\bar t) \bigr) \cap \Big\{\gamma(t)\,;\quad 
 t \in [\bar t -\ve, \,\bar t+\ve]  \Big\} ~= ~\emptyset. 
\eeq

Indeed, 
by a translation and rotation of coordinates, we can assume that $\bar t = 0$, 
$\gamma(0) = 0$ and $\dot \gamma = (1,0)$.   The equation for $\gamma$ can be locally written as
\begin{equation*}
\dot x(t) = \cos \theta(t), \quad \dot y(t) = \sin \theta(t), \qquad - \frac{\pi}{2} \,<\,\theta(0-)\,\leq \,\theta(0+)\,<\, \frac{\pi}{2},
\end{equation*}
where the map $t\mapsto \theta(t)$  
satisfies (\ref{tsc}).

Since  $t\mapsto x(t)$ is Lipschitz, in the interval of invertibility (i.e.~as long as 
$\cos \theta(t) > 0$) we obtain
\begin{equation*}
\frac{\theta(t_2) - \theta(t_1)}{x(t_2) - x(t_1)} ~\geq ~- \frac{t_2-t_1}{r}\cdot \left(\int_{t_1}^{t_2} \cos \theta(\tau) d\tau\right)^{-1}.
\end{equation*}
This implies that $x \mapsto \theta(t(x))$    satisfies
\begin{equation*}
D_x \theta(t(x))\, \geq \,- \frac{1}{r \cos \theta}, \qquad \quad \sin(\theta(t(x))) \,\geq\,- \frac{x}{r},
\end{equation*}
\begin{equation*}
\frac{dy}{dx}\, =\, \tan \theta(t(x))\, \geq\, - \frac{x}{\sqrt{r^2 - x^2}}.
\end{equation*}
Since this argument is valid both for $t>0$ and for $t<0$, one concludes that 
$$y(t)\, \geq\, r - \sqrt{r^2 - (x(t))^2}$$ 
in a nonempty interval  $t\in [-\ve, \ve]$ where $\cos \theta(t) > 0$. 
This proves (\ref{etee}).

%
%
\v
{\bf 2.}  Next, to prove the global property (\ref{eter}), consider any point $\gamma( t_1)\in \partial E$,
and choose $t_2>t_1$ so  that  $\gamma(t_2)$  is the first point where the curve $\gamma$ touches 
again the set $\clos\,B_\rho\bigl(\gamma(t_1) + \rho \bfn(t_1)\bigr)$.
In other words (see Fig.~\ref{f:sc107}, right),
$$t_2~=~\min\,\Big\{t>t_1\,;~~\gamma(t) \in {\rm clos}\,B_\rho\bigl(t_1 + \rho \bfn(t_1)\bigr)
\Big\}.$$
Notice that, by (\ref{etee}), one has $t_2\geq t_1+\ve$, hence the above minimum is well defined.

If   $ t_2 -  t_1 < \rho\pi$, we will derive a contradiction.  As shown in Fig.~\ref{f:sc107}, right,  let $\bar t \in [ t_1, t_2]$ be such that the distance
$
\bigl| \gamma(\bar t) - (\gamma( t_1) + \rho \bfn( t_1)) \bigr|$
 is maximal.
Consider the ball $B_\rho(z)$, where 
$$z~=~ \gamma(\bar t) + \rho \bfn(\bar t),\qquad\qquad    \bfn(\bar t)~=~{\gamma(t_1) + \rho \bfn(t_1) -\gamma(\bar t) \over
\bigl|   \gamma(t_1) + \rho \bfn(t_1) -\gamma(\bar t)\bigr|}\,.$$
Since $t_2 -  t_1 \leq \pi r$, then one of the intervals $[ t_1,\bar t]$, $[\bar t, t_2]$ is shorter than $| t_2 - t_1|/2 \leq \pi r/2$.  
To fix ideas,  assume $\bar t- t_1\leq t_2-\bar t$.
Starting from $\gamma(\bar t)$, we move along this shorter arc  $[ t_1,\bar t]$ toward 
$\gamma(t_1)$ until we reach a first point
$\gamma(\tau )\in \partial B_\rho\bigl(\gamma(\bar t) + \rho \bfn(\bar t)\bigr)$.
Notice that this first intersection point is well defined and bounded away from $\gamma(\bar t)$,  because of the argument in step {\bf 1}. 

Moreover, $\tau $ cannot coincide with $ t_1$.   Indeed, the point $\gamma(\tau)$ lies on the  half circumference
$$\bigl\{y\in B_\rho (z)\,;~(z-y)\cdot \bfn(\bar t) \geq 0\bigr\},$$
while $\gamma(t_1)$ lies on the arc of circumference
$$\Big\{y\in B_\rho \bigl(\gamma(t_1) + \rho \bfn(t_1)\bigr)\,;~(z-y)\cdot \bfn(\bar t) \geq 0\Big\}.$$
These two arcs do not have any point in common.

We then repeat the above construction, replacing $[t_1, t_2]$ with this new interval 
$[\bar t, \tau ]$.
By induction, we thus obtain a sequence of nested intervals $[ t_{1,n}, t_{2,n}]$, 
with 
$ t_{2,n} -  t_{1,n} < 2^{-n} \pi r$, such that
\begin{equation*}
\gamma( t_{2,n}) \in \partial B_\rho(\gamma( t_{1,n}) + \rho \bfn( t_{1,n})), \quad \gamma\bigl(\,] t_{1,n}, t_{2,n}[\,\bigr) \cap 
\clos\,B_\rho(\gamma( t_{1,n}) + \rho \bfn( t_{1,n}))~ = ~\emptyset.
\end{equation*}
Since the length of these intervals is converging to zero, if $t^*$ is the limit of the sequence, then for $n \gg 1$ we obtain a contradiction with 
step {\bf 1},  since $[t_{1,n}, t_{2,n}]$ would be contained in the set where $\gamma$ 
does  not intersect any of the tangent open balls.
This concludes the proof.
\endproof

\begin{corollary} In the above setting, if  (\ref{tsc}) holds and 
${\mathcal H^1}(\partial E) \leq 2\pi r$, then $E$ is 
$r$-semiconvex. More precisely, each point $\gamma(t)+ r\bfn(t)$ has distance $r$ from $E$, and $E = (E^r)^{-r}$.
\end{corollary}

Moreover, for every portion of the  boundary $\gamma$ of length $t_2-t_1\leq \pi r$, the rays $$\rho \mapsto \gamma(t) + \rho \bfn(t)\,,\qquad\qquad \rho\in [0, r]\,,~~t\in [t_1, t_2]$$ 
can intersect only at the initial point $\rho = 0$ and final point $\rho = r$, at most. We can
thus obtain the same formulas  (\ref{32})-(\ref{33}) 
over this portion of $\partial E$.  
 For every Borel subset $A \subset \gamma([t_1, t_2])$, there holds
\begin{equation}
\label{Equa:boundary_r_semi_1}
{\mathcal H^1} \bigg( \bigcup_{x \in A} x + h \bfn(x) \bigg) ~= ~{\mathcal H^1}(A) + h |D\bfn|(A) ~=~ {\mathcal H^1}(A) + h D\theta(A),
\end{equation}
\begin{equation}
\label{Equa:area_r_semi_1}
{\mathcal L}^2 \bigg( \bigcup_{x \in A} x + h \bfn(x) \bigg)~ = ~h {\mathcal H^1}(A) + \frac{1}{2} h^2 |D\bfn|(A) ~=~ h {\mathcal H^1}(A) + \frac{1}{2} h^2 D\theta(A).
\end{equation}
%
%
%
\subsection{Local perturbations of convex sets.}
Let $E$ be a compact convex set, with boundary parameterized by $t\mapsto \gamma(t)$.
As before, let $e^{i\theta(t)}=\bfn(t)= - \dot\gamma(t)^\perp$.
Let $\bar t$ be a Lebesgue point for $D\theta(t)$, so that 
\begin{equation*}\lim_{\delta\to 0+}{1\over\delta} 
\bigl| D\theta - \kappa \mathcal L^1 \bigr| \bigl([ \bar t -\delta,~\bar t+\delta] \bigr) ~=~0.
\end{equation*}
Here $\kappa\geq 0$ is a constant describing the local curvature.
This implies that for every $\ve>0$ one can find $\delta>0$ such that
\begin{equation*}
\big| \bfn(t) - \bfn(\bar t) - \kappa \dot \gamma(\bar t) (t - \bar t) \big| \,\leq \,\ve |t - \bar t| \qquad \text{for} \ t \in [\bar t -\delta,\,\bar t+\delta].
\end{equation*}

%

Observe that, up to a rigid motion and decreasing $\delta$ in case, locally we can write the curve $\gamma$ as the graph of a positive convex function $f(x)$, $x \in [-\delta,\delta]$, with
\begin{equation}
\label{Equa:varphi_prop}
f(0) = 0, \qquad \big| f'(x) - \kappa x \big| \leq \ve |x|, \qquad \lim_{\delta \to 0} \frac{1}{\delta} \big| D^2 f - \kappa \mathcal L^1 \big|([-\delta,\delta]) = 0.
\end{equation}
Consider a semiconvex function $\phi \leq  0$ such that
\begin{equation}
\label{Equa:phi_prop}
\supp\, \phi \subset [-\delta,\delta], \qquad \qquad D^2 \phi \geq - \frac{1}{r} \mathcal L^1.
\end{equation}
It is easy to see that these conditions implies that $|\phi'(x)| \leq \delta$.

Define the set $E'$ by

\begin{equation}\label{E'}
E' = E \cup \big\{y \in f(x) + [\phi(x),0], x \in [-\delta,\delta] \big\}.
\end{equation}

Its boundary is the simple closed curve $\gamma'$ obtained by replacing the part of its graph equal to $\{(x,f(x){\color{blue})}, x \in [-\delta,\delta]\}$ with the graph $\{(x,f(x) + \phi(x)), x \in [-\delta,\delta]\}$. In the interval $x \in [-\delta,\delta]$, the curvature of $\gamma'$ is clearly a measure and its a.c. part can be easily computed as
\begin{equation*}
\frac{f'' + \phi''}{\sqrt{1 + (f' + \phi')^2}} \geq - \frac{1}{r}.
\end{equation*}

%
%

%
%

Next, let $t$ be a parametrization of the new curve $\gamma'$, and for any choice of  $t$ and   $\bfm= e^{i\alpha}$
with $ \alpha\in  [\theta'(t-),\theta'(t+)]$, consider the segment
\begin{equation*}\bigl\{\gamma(t) + \sigma \bfm\,;\quad 0<\sigma<r\bigr\}.
\end{equation*}
We claim that all these segments
are disjoint.  Indeed,  it is enough to verify this statement for segments when $\gamma'(t)$ belongs to the graph of $f + \phi$ and $\phi \not= 0$. In this case, since the length of the arc is $\mathcal O(\delta) < \pi r$, we can apply the analysis in step {\bf 1} of the proof of Lemma~\ref{l:23}.

By \eqref{Equa:area_r_semi_1}, which can be applied to the whole $\partial E'$ because the optimal rays are disjoint, we obtain 
\begin{equation*}
\begin{split}
{\mathcal L^2((E')^r)} - {\mathcal L^2(E^r)} &=~ {\mathcal L^2(E')} - {\mathcal L^2(E)} + r \big( \mathcal{H}^1(\partial E') - {\mathcal H^1}(\partial E) \big) \\
&= ~\int {- \phi(x) \, dx} + r \int {\big( \sqrt{1 + (f'(x) + \phi'(x))^2} - \sqrt{1 + (f'(x))^2} \big) \, dx}.
\end{split}
\end{equation*}
If  $\bar t$ is a Lebesgue point of $D\theta$ where the  derivative is
$\dot \theta(\bar t) =\kappa {~=\frac{1}{\rho}~ }\geq 0$, {corresponding to the Lebesgue point $x=0$ for $D^2 f$ where $D^2 \phi(0) = \phi''(0) = \kappa = \frac{1}{\rho}$,} for $\delta \ll 1$ we have
\bel{fper}
\begin{split}
\bigl({\mathcal L^2((E')^r)} - {\mathcal L^2(E^r)}\bigr) - \bigl({\mathcal L^2(E')} - {\mathcal L^2(E)}\bigr) &=~ r \int {\big( \sqrt{1 + (f'(x) + \phi'(x))^2} - \sqrt{1 + (f'(x))^2} \big) \, dx} \\
&= r \int {\frac{(f'(x) + \phi'(x))^2 - (f'(x))^2}{\sqrt{1 + (f'(x) + \phi'(x))^2} + \sqrt{1 + (f'(x))^2}} \, dx}  \\
&= r {\int_{- \delta}^{\delta} \bigg( \frac{\bigl(\phi'(x)\bigr)^2}{2} + \phi'({\color{blue}x}) \frac{x}{\rho}  \bigg) \bigl(1+ \mathcal O(\delta) \bigr)\, dx}.
\end{split}
\eeq
To achieve the last estimate in (\ref{fper}), notice that {$\phi' = \mathcal O(\delta), \,f' = \frac{x}{\rho} (1 + \mathcal O(\delta))$, because of (\ref{Equa:varphi_prop},\ref{Equa:phi_prop})}.
\v
A useful  choice of the perturbation is
\bel{pper}
\phi({x}) ~=~\left\{\bega{cl} \ds - \frac{\bigl(\delta - |{x}|\bigr)^2}{{2}a}\quad &\hbox{if}~~|{x}|<\delta,
\\[4mm]
0\quad &\hbox{if}~~| x|\geq \delta,
\enda\right.\eeq
with ${a>r}$. In this case, if the local curvature is $\kappa = {1\over\rho}>0$, 
then the perturbed set $E'$ in (\ref{E'}) satisfies
\begin{equation*}
\bega{rl} 
{\mathcal L^2((E')^r \setminus E^r)} - {\mathcal L^2(E' \setminus E)} &=~\ds r {\int_{-\delta}^{\delta} \bigg( \frac{\bigl(\phi'(x)\bigr)^2}{2} + \phi'(t) \frac{x}{\rho}  \bigg) \bigl(1+ \mathcal O(\delta) \bigr)\, dx} \\[4mm]
&=\ds~ r {(1 + \mathcal O(\delta)) \int_{-\delta}^{\delta} \bigg[ \frac{(\delta - |x|)^2}{2a^2} - \frac{|x|(|\delta - |x|)}{a \rho} \bigg]\,dx} \\[4mm]
&=~\ds r \bigl(1 + o(1)\bigr) \bigg[ \frac{\delta^3}{3 a^2} + \frac{\delta^3}{3a\rho} \bigg] ~=~ r 
\bigl(1+o(1)\bigr) \bigg( \frac{1}{a} + \frac{1}{\rho} \bigg) \frac{\delta^3}{3a} \\[4mm]
&= ~\ds r\bigl(1 + o(1)\bigr) \bigg( \frac{1}{a} + \frac{1}{\rho} \bigg) {\mathcal L^2(E' \setminus E)}.
\enda
\end{equation*}
In particular, by letting $a \to +\infty$, we obtain {the following proposition:

\begin{proposition}
\label{Prop:perturb_convex_add}
If $E$ is a convex set such that there is a Lebesgue point for the curvature $\kappa = \frac{1}{\rho}$, then for every $\ve > 0$, there is a set $E' \supset E$ such that
\begin{equation}
\label{E:p1convex}
\bigg| {\mathcal L^2((E')^r \setminus E^r)} - \bigg( 1 + \frac{r}{\rho} \bigg) {\mathcal L^2(E' \setminus E)} \bigg|
~ <~ \ve\, {\mathcal L^2(E' \setminus E)}.  
\end{equation}
\end{proposition}
}

\v

Next, we study what happens when we remove a set $E'$ from a convex set $E$, so that the difference $E \setminus E'$ is still convex. Consider a point $\bar t$ such that
\begin{equation*}
\liminf_{t \to \bar t} ~\frac{\theta(t) - \theta(\bar t)}{t - \bar t} ~\geq ~\frac{1}{\rho}.
\end{equation*}

\begin{figure}[ht]
\centerline{\hbox{\includegraphics[width=17cm]{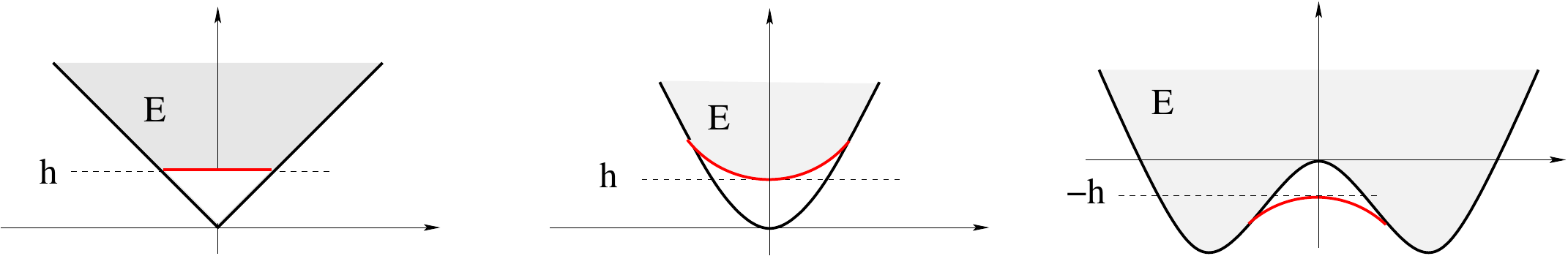}}}
\caption{\small Left and center: the perturbations of the convex set $E$ considered at (\ref{E1h}) and at 
(\ref{E2h}). Right: the perturbation of the semiconvex set $E$ considered at (\ref{Eph}).
  }
\label{f:sc108}
\end{figure}

We first consider the case where $\bar t$ is a corner point. Up to a change of coordinates (see Fig.~\ref{f:sc108}), 
the boundary $\partial E$ is thus the graph of a function $x_2 = f(x_1)$, with
\begin{equation*}
\lim_{x_1 \to 0\pm} f'(x_1) \,=\, \pm a, \qquad a > 0.
\end{equation*}
We then define (see Fig.~\ref{f:sc108}, left)
\bel{E1h}
E'(h)) ~=~ E \cap \{x_2 \geq h\}.
\eeq
Using  \eqref{33}, one obtains
\begin{equation*}
{\mathcal L^2(E \setminus E'(h))} ~=~ \frac{h^2}{a}\bigl(1 + o(1)\bigr),
\end{equation*}
\begin{equation}
\label{222}
\begin{split}
{\mathcal L^2(E^r \setminus (E'(h))^r)} - {\mathcal L^2(E \setminus E'(h))} &~=~ r \big( {\mathcal H^1}(\partial E) - {\mathcal H^1}(\partial E'(h)) \big) \\
&~=~ r \bigl( \sqrt{1 + a^2} - 1 \bigr) \frac{2h}{a} \bigl(1 + o(1)\bigr) \\
&~= ~\frac{2r}{h} \Big( \sqrt{1 + a^2} - 1 \Big) {\mathcal L^2(E \setminus E'(h))} \bigl(1 + o(1)\bigr).
\end{split}
\end{equation}

A similar computation can be done in the case
\begin{equation*}
\lim_{x_1 \to 0} f'(x_1) = 0, \quad \lim_{x_1 \to 0} \frac{f'(x_1)}{x_1} = + \infty. 
\end{equation*}
In this case, {indeed,} for $h \ll 1$ one has $|f'(x_1)| > k |x_1|$, $f(x_1) < h$, with $k \gg 1$ arbitrarily
large. Defining again $E'(h) = E \cap \{x_2 \geq h\}$, one finds
\begin{equation}
\label{224}
\bega{l}
{\mathcal L^2(E^r \setminus (E'(h))^r)} - {\mathcal L^2(E \setminus E'(h))} ~=~ r \Big( {\mathcal H^1}(\partial E) - {\mathcal H^1}(\partial E'(h)) \Big) \\[4mm]
\ds\qquad =~ r \int_{f(x_1) < h} \bigl( \sqrt{1 + (f'(x_1))^2} - 1 \bigr) \,dx_1 
~\geq~ \frac{r}{2} (1 + o(1)) \int_{f(x_1) < h} (f'(x_1))^2 dx_1 \\[4mm]
\qquad\ds \geq \frac{rk}{2} \int_{f(x_1) < h} x_1 f'(x_1) ~=~ \frac{rk}{2} {\mathcal L^2(E \setminus E)}. 
\enda
\end{equation}

Finally, assume  that 
\begin{equation*}
\lim_{x \to 0\pm} \frac{f'(x)}{x} = \frac{1}{\rho},
\end{equation*}
so that $f(x_1) = \frac{x_1^2}{2\rho}\bigl(1 + o(1)\bigr)$. Choosing $a>\rho$ and defining
(see Fig.~\ref{f:sc108}, right)
\bel{E2h}
E'(h) = E \cap \bigg\{x_2 \ge g(x_1) \doteq h + \frac{(x_1)^2}{2a} \bigg\},
\eeq
we now obtain
\begin{equation*}
{\mathcal L^2(E \setminus E'(h))} ~= ~\int_{f < g} x_1 \bigl( f'(x_1) - g'(x_1) \bigr) dx_1 ~
=~ \bigl(1 + o(1)\bigr) \bigg( \frac{1}{\rho} - \frac{1}{a} \bigg) \int_{f < g} (x_1)^2 dx_1, 
\end{equation*}
\begin{equation*}
\begin{split}
{\mathcal L^2(E^r \setminus (E'(h))^r)} - {\mathcal L^2(E \setminus E'(h))} &= r \big( {\mathcal H^1}(\partial E) - {\mathcal H^1}(\partial E'(h)) \big) \\
&= r \int_{f < g} \bigg( \sqrt{1 + \bigl(f'(x_1)\bigr)^2} - \sqrt{1 + \bigl(g'(x_1)
\bigr)^2} ) \bigg) dx_1 \\
&= \frac{r}{2} \bigl(1 + o(1)\bigr)\int_{f < g} \bigl( (f'(x_1))^2 - (g'(x_1))^2 \bigr) dx_1 \\
&= \frac{r}{2} \bigl(1 + o(1)\bigr) \bigg( \frac{1}{\rho^2} - \frac{1}{a^2} \bigg) \int_{f < g} x_1^{2} \, dx_1 \\
&= \frac{r}{2} \bigl(1 + o(1)\bigr) \bigg( \frac{1}{\rho} + \frac{1}{a} \bigg) 
{\mathcal L^2(E \setminus E'(h))}.
\end{split}
\end{equation*}
For every $\ve>0$,
letting $a \nearrow \rho$ we obtain that,  for $h >0$ small enough,
\begin{equation}
\label{Eper2}
\bigg| {\mathcal L^2(E^r \setminus (E'(h))^r)} - \bigg( 1 + \frac{r}{\rho} \bigg) 
{\mathcal L^2(E \setminus E'(h))} \bigg| ~< ~\ve\, {\mathcal L^2(E \setminus E'(h))}.
\end{equation}

{We summarize the results into the following proposition:

\begin{proposition}
\label{Prop:perturb_convex_subtract}
If $E$ is a convex set, and there is a point where the inner radius of curvature is $0$, then for every $\rho > 0$ there is a perturbation $E' \subset E$ such that
\begin{equation}
\label{Equa:substract_corner}
{\mathcal L^2(E^r \setminus (E')^r)} - {\mathcal L^2(E \setminus E')}~ >~
 \frac{r}{\rho} {\mathcal L^2(E \setminus E')}. 
\end{equation}



If at a boundary point $x \in \partial E$ the radius of curvature is $\rho > 0$, then for every $\ve > 0$  there is a set $E' \subset E$ such that
\begin{equation}
\label{E:p1convex_2}
\bigg| {\mathcal L^2}\bigl(E^r \setminus (E')^r\bigr) - \bigg( 1 + \frac{r}{\rho} \bigg) {\mathcal L^2(E\setminus E')} \bigg|
~ <~ \ve\, {\mathcal L^2(E \setminus E')}.  
\end{equation}
\end{proposition}
}

\v
\subsection{Local perturbations of $r$-semiconvex sets.}

Similar computations can be done locally for a semiconvex set $E$. To fix ideas, consider a boundary point $\bar x = \gamma(\bar t)\in \partial E$. As before, denote by 
$e^{i\theta(t)} =\bfn(t)$ the set of outer normals and assume that  $\bar t$ is a Lebesgue point of $\theta$, with $\dot \theta (\bar t)= - \frac{1}{\rho}$.  
Writing $\partial E$ locally as the graph of a semiconvex function $x_2 = f(x_1)$ with $\ddot f(0) = - 1/\rho$, we choose $a>\rho$ and  replace $E$ with the slightly larger sets
\bel{Eph}
E'(h) ~\doteq~
~E \cup \bigg\{x_2 \geq  g(x_1)\doteq- h - \frac{x_1^2}{2 a} \bigg\}.
\eeq

Similarly to the previous cases,   we compute
\begin{equation*}
\begin{split}
{\mathcal L^2((E'(h))^r \setminus E^r)} - {\mathcal L^2(E'(h) \setminus E)} &~\leq 
~r \bigl( {\mathcal H^1}(\partial E'(h)) - {\mathcal H^1}(\partial E) \bigr) \\
&~=~ - r \int_{f < g} \bigg( \sqrt{1 + (f'(x_1))^2} - \sqrt{1 + (g'(x_1))^2} ) \bigg) dx_1 \\
&~=~ - \frac{r}{2} \bigl(1 + o(1)\bigr) \int_{f < g} \big( (f'(x_1))^2 - (g'(x_1))^2 \big) dx_1 \\
&~=~ - \frac{r}{2} \bigl(1 + o(1)\bigr) \bigg( \frac{1}{\rho^2} - \frac{1}{a^2} \bigg) \int_{f < g} (x_1)^2 dx_1 \\
&~=~ - \frac{r}{2} \bigl(1 + o(1)\bigr) \bigg( \frac{1}{\rho} + \frac{1}{a} \bigg) 
{\mathcal L^2(E'(h)\setminus E)}.
\end{split}
\end{equation*}
The first inequality is due to the fact that \eqref{Equa:area_r_semi_1} holds only locally, and in general there can be point in $E^r$ belonging to more than one optimal ray. Letting $a \nearrow \rho$ we obtain {the following lemma:

\begin{proposition}
\label{Lem:perturb_semiconvex}
Let $E$ be a semiconcave set and $x$ a point with outer curvature $\rho$. Then for every $\ve>0$ there exists a set $E' \supset E$ such that 
\begin{equation}
\label{Eper3}
\bigg| {\mathcal L^2( (E')^r \setminus E^r)} - \bigg( 1 - \frac{r}{\rho} \bigg) {\mathcal L^2(E' \setminus E)} \bigg| 
~<~ \ve {\mathcal L^2(E' \setminus E)}.
\end{equation}
\end{proposition}
}

\begin{figure}[ht]
\centerline{\hbox{\includegraphics[width=10cm]{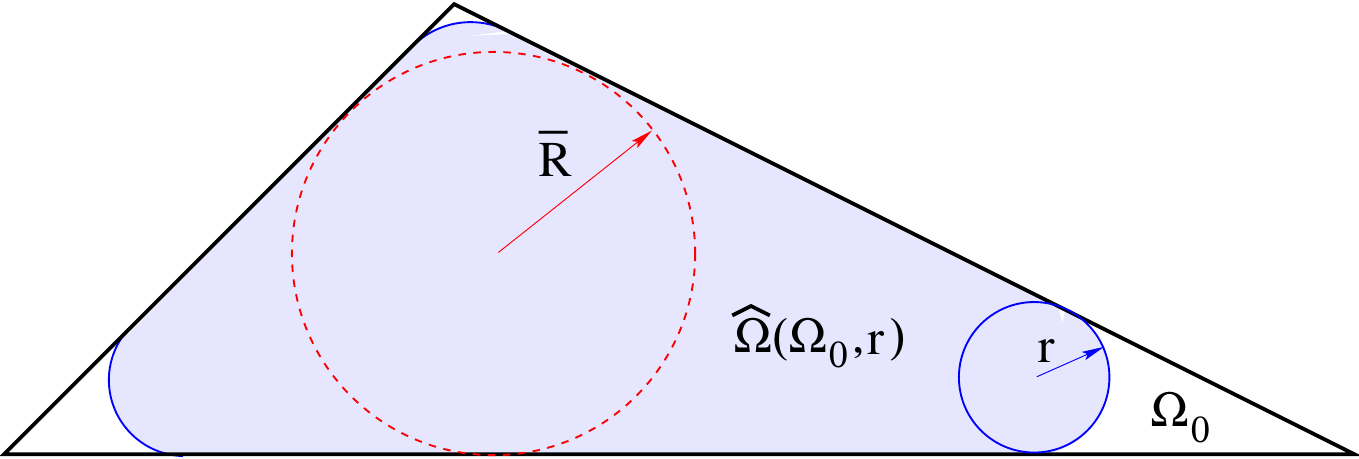}}}
\caption{\small  The maximum inner radius $\ov R$ defined at (\ref{ovrdef}),  and the set $\Hat \Omega(\Omega_0 , \rho)$
(shaded region), in the case where
$\Omega_0$ is a triangle. }
\label{f:sc73}
\end{figure}

\section{A  one-step minimization problem}
\label{s:3}
\setcounter{equation}{0}
Given a bounded convex {closed} set $\Omega_0  \subset \R^2$ with nonempty interior, 
its inner radius is defined by setting
\bel{ovrdef}
\ov  R ~= ~\ov R(\Omega_0)~\doteq~ \max \bigl\{r>0\,;~~ \Omega_0  ~\hbox{contains an open ball $B_r(x)$ of radius $r$} \bigr\}.
\eeq
As shown in Fig.~\ref{f:sc73}, for  $\rho \in \,]0,\ov  R]$ we  also consider the open set
\bel{HOdef}
\Hat \Omega(\Omega_0 ,\rho) ~\doteq~ \bigcup_{B_\rho(x) \subseteq \Omega_0}
 B_\rho(x)\,.
\eeq
In this section, 
given a constant $0<a\leq 
\mathcal L^2(\Omega_0)$, and a radius $r>0$, we will study  the following one-step minimization problem:
\begin{equation}
\label{op1} \hbox{minimize:}~~~ {\mathcal L^2(\Omega^r)},\qquad\hbox{subject to}\quad\Omega 
\subseteq \Omega_0 ,
\quad  {\mathcal L}^2 (\Omega) = a.\end{equation}
In other words, among all sets of fixed area $a$ contained inside $\Omega_0 $, we seek one that minimizes
the area of its $r$-neighborhood.  We will show that the optimal solutions to (\ref{op1}) do 
not depend on the 
radius $r>0$.  Namely, a set $\Tilde\Omega$ is optimal if and only if it solves the
corresponding  minimization problem 
for the perimeter:
\bel{op2} \hbox{minimize:}~~~ {\mathcal H^1({\partial}\Omega)},\qquad\hbox{subject to}\quad\Omega \subseteq \Omega_0 ,
\quad  {\mathcal L}^2 (\Omega) = a.\eeq 
The existence and various properties of optimal solutions to (\ref{op2}) have been established 
in  \cite{SZ}.   To show the equivalence of the two problems (\ref{op1}) and (\ref{op2}),
the key step is to prove that the optimal solutions of (\ref{op1}) are convex.
This is the content of the following theorem.
\v

\begin{theorem}
\label{t:31}
Consider a compact convex set   ${\Omega_0}\subset\R^2$  and let $0 < a \leq |{\Omega_0}|$.
Then there exists a set
\bel{argmin}
\Omega ~=~ \argmin \bigl\{{\mathcal L^2(\Omega^r)}\,;~~\Omega ~\hbox{is a closed  subset of ${\Omega_0}$ with area }~{\mathcal L}^2 (\Omega) = a \bigr\}.
\eeq
Moreover, every such minimizer is convex.
\end{theorem}

The proof will be achieved in several steps. The existence of a minimizer follows from 
a standard compactness argument. However, its convexity 
requires a careful analysis.  First we prove that every connected component of an optimal set $\Omega$ must be convex.  Then we show that an optimal set 
can have at most finitely many components. Finally, we will prove that every optimal set
$\Omega$ is connected.

\subsection{Existence of a minimizer.} 
We prove here that the optimization problem (\ref{argmin}) has a solution.
Let $(\Omega_n)_{n\geq 1}$ be a minimizing sequence of compact sets, such that
\bel{Omn}
\Omega_n\ \subseteq {\Omega_0}, \qquad {\mathcal L^2(\Omega_n)} \geq a, \qquad 
\lim_{n \to \infty} {\mathcal L^2(\Omega_n^r)} ~=~ m~\doteq
\inf \bigl\{{\mathcal L^2(\Omega^r)}\,;~~ \Omega \subseteq {\Omega_0},~ {\mathcal L}^2 (\Omega) = a \bigr\}.
\eeq
By possibly replacing $\Omega_n$ with the larger set 
$(\Omega_n^r)^{-r}$, which is still contained inside ${\Omega_0}$, we can assume that 
$\Omega_n^r$ and $\Omega_n$ are in duality.
%

{We use the following lemma.

\begin{lemma}
\label{l:31H}
If $\Omega_n,\Omega_n^r$, $n \in \N$, are a family of sets in duality, with $\Omega_n \subset \clos\,B_R(0)$ for a fixed closed ball $\clos\,B_R(0)$, then there exists a subsequence $\Omega_{n_k},\Omega_{n_k}^r$, $k \in \N$, and a set $\Omega$ such that
\begin{equation*}
\Omega_{n_k} \to \Omega, \qquad \qquad \clos \, B_{R+r}(0) \setminus \Omega_{n_k}^r \to \clos \, B_{R+r}(0) \setminus \Omega^r, 
\end{equation*}
w.r.t.~the Hausdorff distance of compact sets and w.r.t.~the $\L^1$-distance, respectively.\\
In particular, the limit sets are in duality.
\end{lemma}

{\bf Proof.} {\bf 1.}   By possibly taking a subsequence, we obtain the 
convergence
\begin{equation*}
\Omega_{n_k} \underset{\text{Hausdorff}}{\longrightarrow} \Omega_{\mathrm{Haus}}, \qquad \qquad\Omega_{n_k} \underset{\L^1}{\longrightarrow} \Omega_{L^1},
\end{equation*}
\begin{equation*}
\clos\,B_{R+r}(0) \setminus \Omega^r_{n_k} \underset{\text{Hausdorff}}{\longrightarrow} \clos\,B_{R+r}(0) \setminus \Omega'_{\mathrm{Haus}}, \qquad \qquad\Omega^r_{n_k} \underset{\L^1}{\longrightarrow} \Omega'_{L^1},
\end{equation*}
for some limit sets
$\Omega_{\mathrm{Haus}}$, $\Omega_{L^1}$, $\Omega'_{\mathrm{Haus}}$, and $\Omega'_{L^1}$, all contained in  the closed ball $\clos\,B_{R+r}(0)$. Indeed, this
follows immediately from the finite perimeter estimate (Lemma \ref{l:22}) and the compactness of the family of compact subsets of the compact set $\clos\, B_{R+r}(0)$.
\v
{\bf 2.}
It remains to prove that
$$
\Omega'_{\mathrm{Haus}} = \Omega_{\mathrm{Haus}}^r, \qquad \Omega_{\mathrm{Haus}} = \Omega_{L^1}, \qquad \Omega_{\mathrm{Haus}}^r = \Omega'_{L^1},
$$
where the last two equalities should be intended as $\L^1$-equivalence of the characteristic functions.

First of all, we have that if $y \in \Omega_{\mathrm{Haus}}^r$, then there exists a point $x \in \Omega_{\mathrm{Haus}}$ such that $|x-y| = r - \delta r < r$ for some $\delta r > 0$. Then by Hausdorff convergence, there exists $K \gg 1$ such that for all $k \geq K$ it holds $d(\Omega_{n_k},\Omega_{\mathrm{Haus}}) < \delta r/2$. Then if $x_{n_k} \in \Omega_{n_k} \cap B_{\delta r/2}(x)$ it holds
\[
|y - x_{n_k}| \leq |y-x| + |x- x_{n_k}| < r - \frac{\delta r}{2}.
\]
We thus conclude that $y \in \Omega_{n_k}^r$ for all $k$ large enough. 
Hence $\Omega^r_{\mathrm{Haus}} \subset \Omega'_{\mathrm{Haus}}$. 

A similar argument shows that, if $y \notin \clos\, \Omega_{\mathrm{Haus}}^r$, then $y \notin \Omega_{n_k}^r$ for all $k$ large enough. 
Moreover, for every $\ve > 0$ we obtain
\begin{equation}
\label{Equa:limit_Haus}
(\Omega_{\mathrm{Haus}}^r)^{-\ve} \subset \Omega_{n_k}^r \subset \Omega_{\mathrm{Haus}}^{r+\ve}
\end{equation}
for all $k$ suitably large. Since $\Omega_{\mathrm{Haus}}^r$ is open, we conclude that
\begin{equation*}
\Omega'_{\mathrm{Haus}} = \Omega_{\mathrm{Haus}}^r.
\end{equation*}

Next, since the boundary of $\Omega_{\mathrm{Haus}}^r$ is rectifiable, we have 
\begin{equation*}
\mathcal L^2 \bigl( \Omega_{\mathrm{Haus}}^{r+\ve} \setminus (\Omega_{\mathrm{Haus}}^r)^{-\ve} \bigr) ~=~ \mathcal L^2 \bigl( \partial \Omega_{\mathrm{Haus}}^r + B_\ve(0) \bigr)~ = ~
(2 \ve + o(\ve)) \mathcal H^1(\Omega_{\mathrm{Haus}}^r),
\end{equation*}
and from \eqref{Equa:limit_Haus} we get that the symmetric difference satisfies
\begin{equation*}
\lim_{k\to\infty} \mathcal L^2 \big( \Omega_{n_k}^r \Delta \Omega_{\mathrm{Haus}}^r \big) = 0,
\end{equation*}
i.e. $\Omega_{L^1}' = \Omega_{\mathrm{Haus}}^r$.

Reversing the analysis, i.e. considering $\Omega_{n_k} = (\Omega_{n_k}^r)^{-r}$, we obtain that for every $\ve > 0$ there exists $K > 0$ such that, for $k \geq K$, 
\begin{equation*}
(\Omega_{\mathrm{Haus}}^r)^{-r-\ve} ~\subset ~(\Omega_{n_k}^r)^{-r}~ \subset ~(\Omega_{\mathrm{Haus}}^r)^{-r+\ve}.
\end{equation*}
Being $(\Omega_{n_k}^r)^{-r} = \Omega_{n_k}$, we obtain up to negligible sets
\begin{equation*}
\inter\,\Omega_{\mathrm{Haus}} ~\subset ~\Omega_{L^1} ~\subset~ (\Omega_{\mathrm{Haus}}^r)^{-r} 
~=~ \Omega_{\mathrm{Haus}}.
\end{equation*}
This yields  $\Omega_{L^1} = \Omega_{\mathrm{Haus}}$, because $\mathcal L^2(\partial \Omega_{\mathrm{Haus}}) = 0$.
\endproof

}

%
%
Applying this result to the sequence  introduced at (\ref{Omn}), we obtain a subset $\Omega = \lim_{k} \Omega_{n_k}$ such that $\Omega^r = \lim_{k} \Omega_{n_k}^r$ and
$${\mathcal L}^2 (\Omega)\,=\,
\lim_{{k}\to\infty} {\mathcal L^2(\Omega_{n_k})}\,\geq\, a,\qquad \qquad {\mathcal L^2(\Omega^r)}\,=\,\lim_{{k}\to\infty} {\mathcal L^2(\Omega_{n_k}^r)}\,=\,m.$$

It remains to show that ${\mathcal L}^2 (\Omega)=a$.  If on the contrary ${\mathcal L}^2 (\Omega)>a$, we fix a unit vector $\bfe_1\in \R^2$
and consider the smaller sets
\bel{Olam}
\Omega(\lambda)~\doteq~\bigl\{x\in\Omega\,;~~x\cdot \bfe_1\,\leq\,\lambda\bigr\}.\eeq
For a suitable $\lambda\in \R$, one has ${\mathcal L^2(\Omega(\lambda))}= a$.
But this would imply ${\mathcal L^2(\Omega(\lambda)^r)} < {\mathcal L^2(\Omega^r)} = m$,
reaching a contradiction.

{We collect the results of this section into the following proposition.

\begin{proposition}
\label{Prop:minimal_optimal}
Every minimizer $\Omega$ of (\ref{argmin}) satisfies
\begin{equation*}
\Omega = (\Omega^r)^{-r} \quad \text{and} \qquad \Omega = \clos(\inter\, \Omega).
\end{equation*}
\end{proposition}

{\bf Proof.} The first identity was already proved in Lemma~\ref{l:31H}.  The second identity
can be proved as follows. 
Since $\mathcal L^2(\partial \Omega) = 0$, then $ \mathcal L^2\big(\clos(\inter\, \Omega)\big)= a$. Moreover, for all $x \in \partial (\inter\,\Omega)$ there is $y \in \partial \Omega^r$ such that $|x-y| = r$, being $\partial (\inter\,\Omega) \subset \partial \Omega$.   This yields 
$$
\clos(\inter\,\Omega) = ((\clos(\inter\,\Omega))^r)^{-r}.
$$
If $\Omega \varsupsetneq \clos\,(\inter\, \Omega)$, then $\Omega^r \varsupsetneq (\clos(\inter\,\Omega))^r$.   Since these two sets are open, we would have $\mathcal L^2(\Omega^r) \geq \mathcal L^2((\clos(\inter\,\Omega))^r)$, contradicting the optimality of $\Omega_r$.
\endproof

}

\begin{figure}[ht]
\centerline{\hbox{\includegraphics[width=9cm]{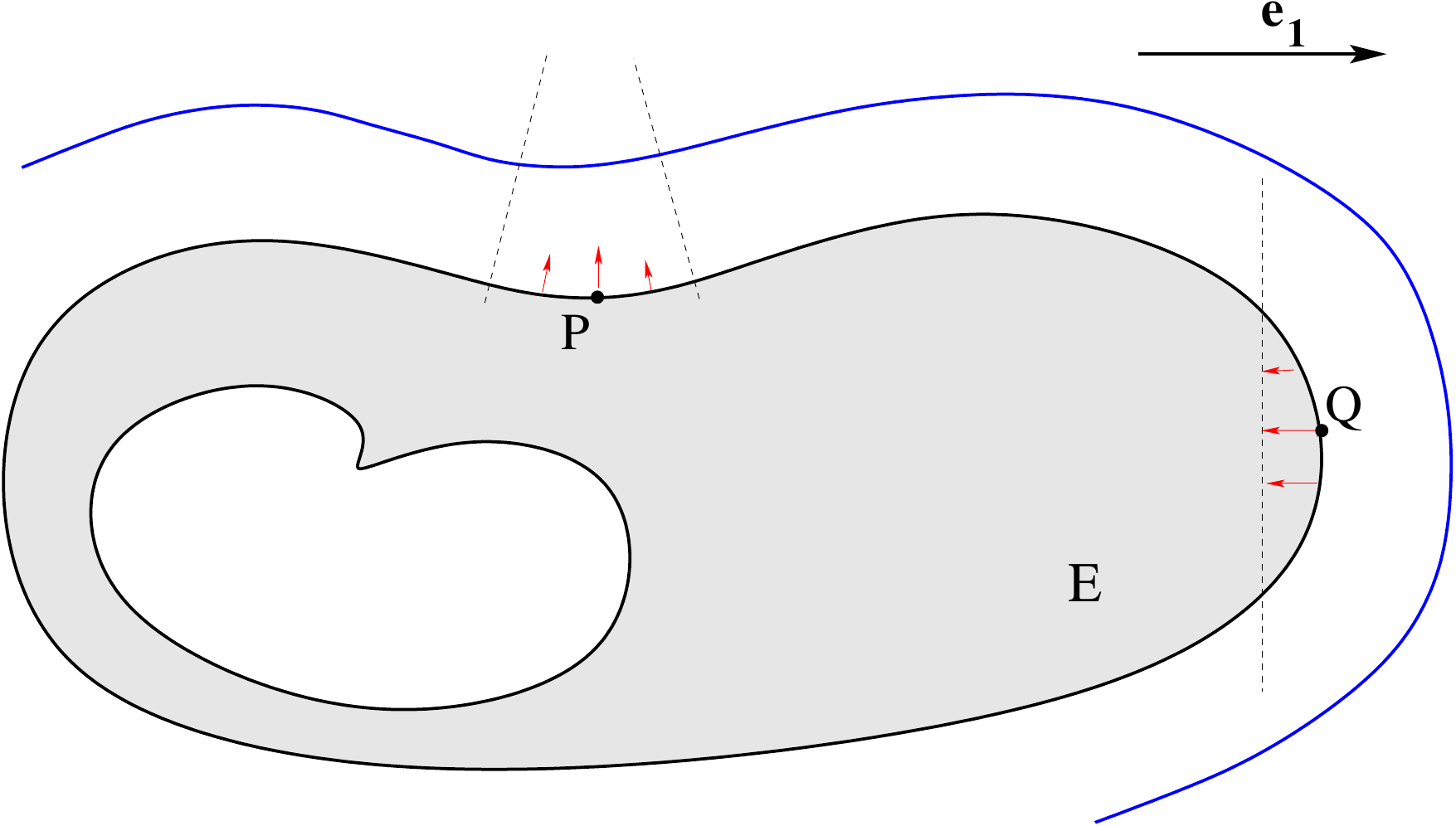}}}
\caption{\small Proving Lemma~\ref{l:32}. If the set $E$ is not convex, we can enlarge it 
in a neighborhood of the point $x^*$ where the curvature of the boundary is negative. 
At the same time, we shrink it in a neighborhood of the exposed point $y$. 
This yields a perturbed set $E_\ve$ with same area as $E$, but with ${\mathcal L^2(E^r_\ve)}<{\mathcal L^2(E^r)}$.  We remark that, in general, the points $y$ and $x^*$ may belong to distinct
connected components of $E$.  }  
\label{f:sc89}
\end{figure}

\subsection{Convexity of the optimal set.}
Aim of this section is to prove that every connected component of a minimizer is a 
compact convex set. Because of Proposition \ref{Prop:minimal_optimal}, it is enough to study the connected components with positive measure. 

\begin{lemma}
\label{Lem:cut_line} Let $\Omega$ be any compact set and let 
$\Omega(\lambda )$ be as  in (\ref{Olam}). 
If ${\mathcal L^2(\Omega(\lambda))} <{\mathcal L^2(\Omega)}$, then
\bel{sdc}
{\mathcal L^2(\Omega^r \setminus \Omega(\lambda)^r)} ~\geq~{\mathcal L^2(\Omega \setminus \Omega(\lambda))}.
\eeq
\end{lemma}

{\bf Proof.}  The inequality (\ref{sdc}) is an immediate consequence of the inclusions
\begin{equation*}
\bigl(\Omega \setminus \Omega(\lambda )\bigr) + r \bfe_1~ \subseteq~ \Omega^r \cap \{x\,;~x\cdot \bfe_1 \geq \lambda  + r\} ~\subseteq ~\Omega^r \setminus \Omega(\lambda )^r.
\end{equation*}\endproof
\v
\begin{lemma} \label{l:32} Let $\Omega\subseteq {\Omega_0}$ be a minimizer for (\ref{argmin}).
Then every connected component of ${\inter\,\Omega}$ is convex.
\end{lemma}


{\bf Proof.} Assume, on the contrary, that
the optimal set $\Omega$ has a connected component  $E$ which is not convex.
Since $\Omega$ is $r$-semiconvex,
according to Lemma~\ref{l:31} there is a point $x^*\in \partial E$ where the boundary
has negative curvature.   We will derive a contradiction showing that $\Omega$ is not optimal.
As shown in  Fig.~\ref{f:sc89}, by slightly enlarging the set $E$ near $x^*$ and shrinking the set $\Omega$ at some other point $y$
along its boundary, we can keep constant the area ${\mathcal L}^2 (\Omega)$, but decrease the area ${\mathcal L^2(\Omega^r)}$ of the $r$-neighborhood.
\v
{\bf 1.}
To construct these perturbations, 
as stated in Lemma~\ref{l:31} we can find  a boundary point $x^*= \gamma_k(t^*)$
where the derivative  $D\theta_k(t^*)$ exists and is strictly negative, say 
$$D\theta_k(t^*) ~=~- {1\over  \rho}\,,$$
with $\rho>0$. 
Constructing  the slightly larger sets $E'(h)$ as in (\ref{Eph}), the change in the area of 
the $r$-neighborhoods 
$(E'(h))^r$ is bounded above by (\ref{Eper3}).
\v
{\bf 2.} Next, for any ${h} > 0$ small enough, choose $\lambda = \lambda_{h}$ so  that 
the set $\Omega(\lambda)$ in (\ref{Olam}) satisfies
\bel{ole}{\mathcal L^2(\Omega \setminus \Omega(\lambda_h))}~=~{\mathcal L^2(E'(h) \setminus E)}.\eeq
Then define the perturbed set
\bel{Ome}
\Omega_{h}~\doteq~(\Omega \cup {E'(h)})\cap \{x\,;~x\cdot \bfe_1\leq \lambda_{h}\}.\eeq
The above definition implies
${\mathcal L^2(\Omega_h)}={\mathcal L}^2 (\Omega)$ for every ${h}>0$ small enough.   Moreover, combining
(\ref{Eper2}) with (\ref{sdc}) we obtain
$$\bega{rl} {\mathcal L^2(\Omega_h^r)}-{\mathcal L^2(\Omega^r)}&\ds\leq~{\mathcal L^2((E'(h))^r \setminus E^r)} -
{\mathcal L^2(\Omega^r \setminus \Omega(\Lambda_h)^r)} \\[4mm]
&\ds \leq ~
 \bigg( 1 - \frac{r}{\rho} + o(1) \bigg) {\mathcal L^2(E(h) \setminus E)} - {\mathcal L^2(\Omega \setminus \Omega(\Lambda_h))}\\[4mm]
 &\ds =~ \bigg( 1 - \frac{r}{\rho} + o(1) \bigg) {\mathcal L^2(E(h) \setminus E)} -  {\mathcal L^2(E(h) \setminus E)} \\[4mm]
&\ds =~ \bigg( - \frac{r}{\rho} + o(1) \bigg) {\frac{h^3}{3}}
 ~<~0\enda
 $$
for all ${h}>0$ small enough. This contradicts the optimality of $\Omega$, proving the lemma.
\endproof

{We observe that the same proof can be adapted to the case of a connected component with $0$ measure. Indeed in this case one observes that every connected set of finite length is covered by a closed curve. We will not need this fact, because we will prove in Lemma \ref{l:36} that there are only finitely many components.}

\begin{figure}[ht]
\centerline{\hbox{\includegraphics[width=12cm]{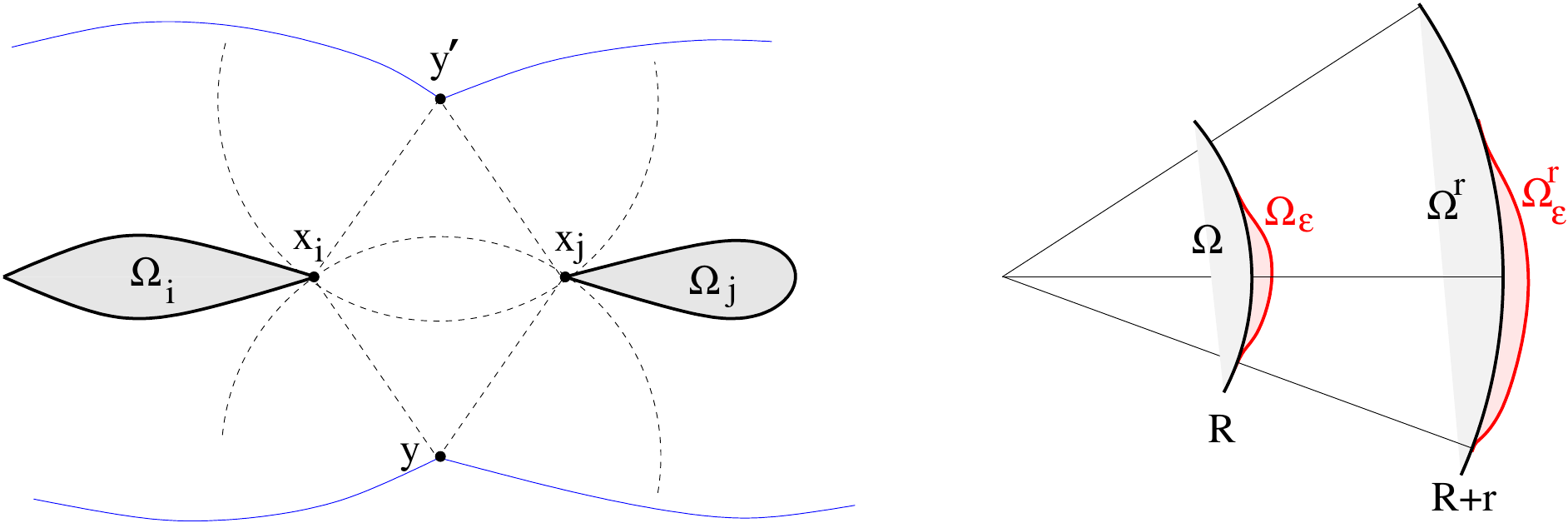}}}
\caption{\small Left: the points $x_i, x_j, y, y'$ considered in Lemma~\ref{l:33}. 
Right: A slight perturbation of the set $\Omega$. If $R$ is the constant curvature radius, then 
the increase in the area of the neighborhood $\Omega_r$ satisfies
${\mathcal L^2(\Omega^r_\ve)}-{\mathcal L^2(\Omega^r)}\approx {R+r\over R} \bigl({\mathcal L}^2(  \Omega_\ve)-{\mathcal L}^2 (\Omega)\bigr)$.
}
\label{f:sc92}
\end{figure}

To prove Theorem~\ref{t:31} it remains to prove that the optimal set $\Omega$  is connected.
As an intermediate step, we will 
show that $\Omega$ has at most finitely many connected components.

Let $\{\Omega_i\}_{i\in \N}$ be the connected components of $\inter\,\Omega$. 

We will use the following lemmas, whose proofs are elementary.

\begin{lemma}
\label{l:33}  Let  $\Omega$ be an optimal set, and let $\Omega_i, \Omega_j$ 
be distinct connected components of ${\inter\,\Omega}$.
If ${(\clos\,\Omega_i)^r \cap (\clos\,\Omega_j)^r} ~\not= ~\emptyset$, then there exists points $ y,y'$, such that
\begin{equation*}
\partial \Omega_i^r \cap \partial \Omega_j^r ~=~ \{y,y'\}.
\end{equation*}
Moreover, there exists points $x_i\in \Omega_i$, $x_j\in \Omega_j$ such that
\bel{xpx}
{(\clos\,\Omega_i)^r \cap (\clos\,\Omega_j)^r} ~=~ B_r(x_i) \cap B_r(x_j), \qquad \quad \{x_i, x_j\} ~= ~\partial B_r(y) \cap \partial B_r(y').
\eeq
{Finally, for every $x \in \R^2$
\begin{equation}
\label{Equa:triple_intersection}
\sharp \big\{i : x \in \Omega_i^r \big\} \leq 2.
\end{equation}}
\end{lemma}

{\bf Proof.}
It is clear that the boundaries of the two open convex sets $\partial \Omega_i^r, \partial\Omega_j^r$
intersect exactly at two points $y,y'$ (see Fig.~\ref{f:sc92}, left).
Let $x_i,  x'_i,\in \partial \Omega_i$ and $ x_j, x'_j\in \partial \Omega_j$ be points such that 
$$|y-x_i| ~=~|y-x_j|~=~|y'-x'_i|~=~|y'-x'_j|~=~r.$$
We claim that $x_i= x'_i$ and $x_j= x'_j$. Indeed, if $x_i\not= x_i'$,
consider the arc along the boundary of 
$\Omega_i$ with endpoints $x_i, x_i'$.  This arc has positive length and thus 
contains at least one point $x\in\partial \Omega_i$ 
where the unit outer normal $\bfn(x)$
is unique.   
By optimality, the point $x+r\bfn(x)$ cannot lie inside the open set $\Omega^r$, otherwise we could enlarge 
the set $\Omega$ in a neighborhood of $x$, without changing $\Omega^r$.
On the other hand, the above construction implies
$x+r\bfn(x)\in \Omega_j^r$, yielding a contradictions.    
A similar argument yields $x_j=x_j'$.

The remaining identities (\ref{xpx}) are clear.

\begin{figure}[ht]
\centering{\resizebox{8cm}{!}{\includegraphics{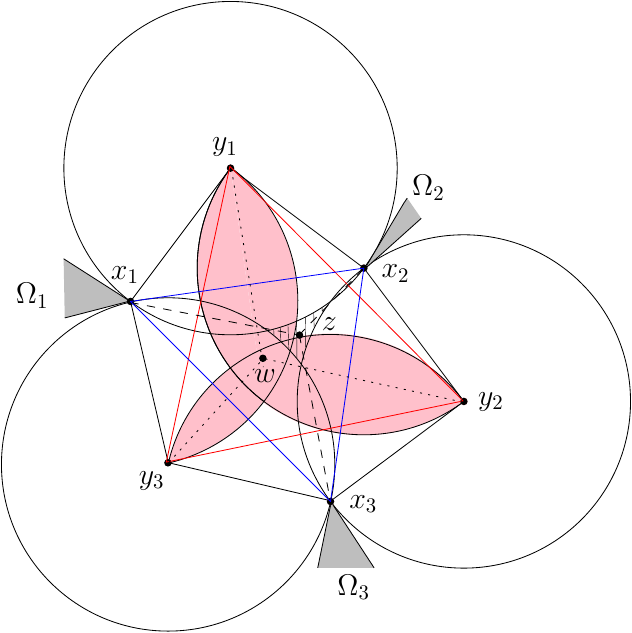}}}
\caption{\small An illustration of the last statement in Lemma \ref{l:33}. If we assume that $w \in \cap_{i=1,2,3} B_r(x_i)$, then the point $z$, intersection of the axis of symmetry of the sides of the triangle $\{y_1,y_2,y_3\}$ cannot belong to any of the balls $B_r(y_i)$, $i=1,2,3$.}
\label{Fig:triple_intersection}
\end{figure}

{
The last estimate follows by the following consideration of elementary geometry. If (\ref{Equa:triple_intersection}) is false, then there is a point $x$ such that (up to relabeling)
\begin{equation*}
x \in \Omega_1^r \cap \Omega_2^r \cap \Omega_3^r,
\end{equation*}
then we are in the situation of Fig. \ref{Fig:triple_intersection}: in particular we can assume that $x$ is the intersection $w$ of the symmetry axis of the sides of the triangle 
$\{x_1,x_2,x_3\}$. The hexagon $\{x_1,y_2,x_2,y_2,x_3,y_2\}$ must be convex, and then
its angles satisfy
\begin{equation}
\label{Equa:exagoon}
\sum_{i=1,2,3} \angle (y_{i-1} x_i y_i) + \sum_{i=1,2,3} \angle (x_i y_i x_{x+1}) = 4\pi, \qquad x_4 = x_1, y_0 = y_3.
\end{equation}
Here and in the following, by $\angle (xyz) $ we denote the angle formed at $y$
by the two segments $xy$ and $yz$.
Since $w \in \cap_{i=1,2,3} B_r(x_i)$, we obtain
\begin{equation*}
\angle (y_{i-1} w y_{i}) ~>~ \angle(y_{i-1} x_{i} y_{i}), \qquad y_0 = y_3,
\end{equation*}
and then
\begin{equation*}
\sum_{i=1,2,3} \angle(y_{i-1} x_i y_{i}) ~<~ \sum_{i=1,2,3} \angle (x_i w x_{i+1}) = 2\pi, \qquad y_0 = y_3.
\end{equation*}
In a similar way, let $z$ be the intersection of the axis of symmetry of the triangle $\{y_1,y_2,y_3\}$: being $z$ equidistant from $y_i$, $i=1,2,3$, we deduce that if $|z - y_i| < r$ then
\begin{equation*}
\sum_{i=1,2,3} \angle(x_i y_i x_{i+1}) ~<~
 \sum_{i=1,2,3} \angle (x_i w x_{i+1}) = 2\pi, \qquad x_4 = x_1.
\end{equation*}
This however contradicts (\ref{Equa:exagoon}). Hence the set $\Omega$ cannot be in duality or optimal: indeed, if it is in duality, the set containing $z$ and not covered by the three balls $B_r(y_i)$, $i=1,2,3$, is not convex, and the only points belonging to $\partial \Omega^r$ at distance $r$ from this set are the points $y_1,y_2,y_3$. 
}
\endproof

In the following, given a distance $\delta>0$, we shall say that two components $\Omega_i,\Omega_j$ are $\delta$-{\bf related} if there exists a point $y$ such that 
\bel{rela}d(y, \Omega_i) ~=~ d(y,\Omega_j)~ =~ d(y,\Omega)~<~\delta.\eeq
By the previous analysis, this can happen only if there are points
$x\in \partial \Omega_i$, $x'\in \partial \Omega_j$ such that, calling $\theta = \bigl|\bfn(x_i)\bigr|$, one has
\bel{cost}r\cos{\theta\over 2} ~<~\delta.\eeq

\begin{corollary}  
For each connected component $\Omega_i$ there are at most 2 sets that are  $(r/2)$-related to $\Omega_i$,
and  at most 3 sets $\Omega_j$ 
which are $(r/\sqrt 2)$-related to $\Omega_i$.
\end{corollary}

{\bf Proof.} The boundary  $\partial \Omega_i$ can contain at most two points where 
$\bigl|\bfn(x)\bigr| > 2\pi/3$.  Taking $\delta =  r \cos{\pi\over 3}={r\over 2}$ by (\ref{cost}) we obtain the first assertion.

Similarly, there can be at most 3 points where $\bigl|\bfn(x)\bigr| > \pi/2$.
Taking $\delta =  r \cos{\pi\over 4}= {r\over \sqrt 2}$,  by (\ref{cost}) we obtain the second assertion.
\endproof

\begin{lemma}\label{l:340} Let $\Omega_i$ be a connected component of {$\inter\, \Omega$, where $\Omega$ is an optimal set}. 
Let
$$\theta_{max}~\doteq~\max_{x\in \partial \Omega_i}\bigl|\bfn(x)\bigr| {< \pi}$$
be the maximal angle at corner points of $\partial \Omega_i$.
Then all other components of $\Omega$ have a strictly positive distance from $\Omega_i$. Namely
\bel{domi}
{(\clos\,\Omega_i)^\rho \cap} \bigl(\Omega\setminus\Omega_i)~=~\emptyset, \qquad\qquad\hbox{with}\quad \rho~=~r\sin {\bigg(}{\pi-\theta_{max}\over 2}{\bigg)}\,.\eeq
\end{lemma}

\begin{figure}[ht]
\centerline{\hbox{\includegraphics[width=13cm]{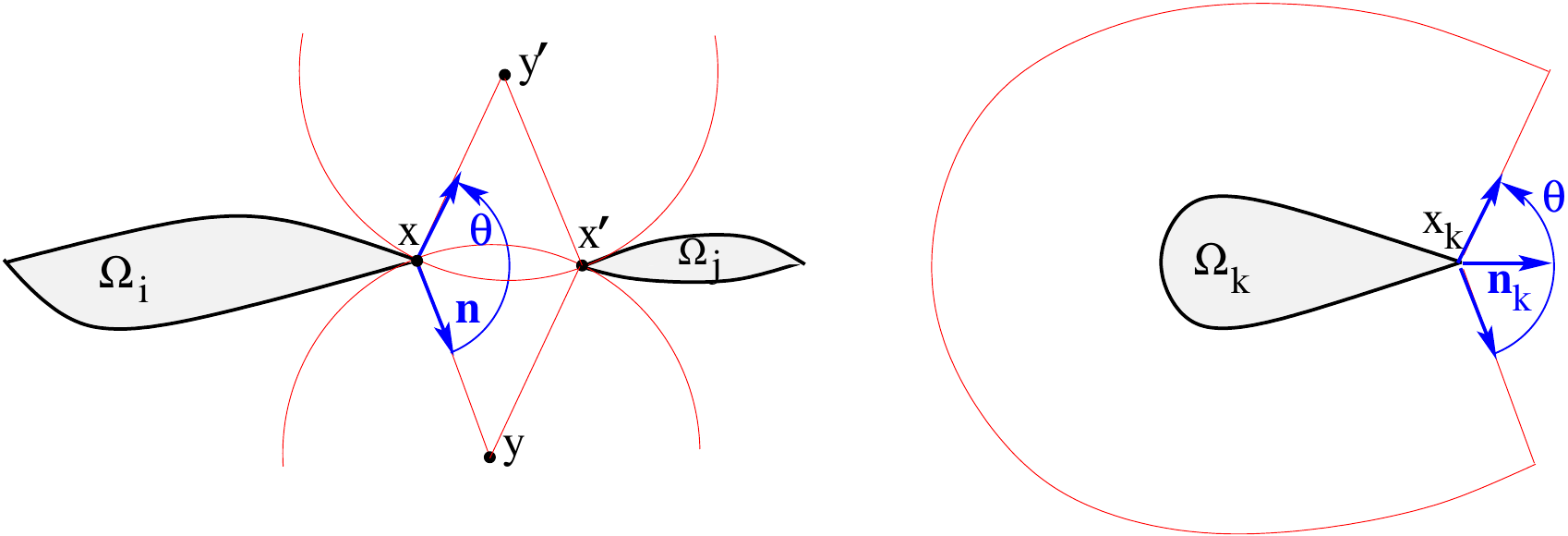}}}
\caption{\small Left: The minimum distance between the component $\Omega_i$ 
and any other component $\Omega_j$
is bounded below in terms of the angle $\theta$. Right: the central unit normal vector $\bfn_k$ at the point $x_k$.
If the set ${\partial \Omega_k}$ does not contain a second point $x_k'$ where the set of outer normals has size 
$\bigl|\bfn(x_k')\bigr|~\approx \pi$, then there is a large  region to the left of $\Omega_k$ which cannot intersect 
any other connected component.}
\label{f:sc101}
\end{figure}

{\bf Proof.} 
Consider the set of boundary points admitting a single outer normal:
\bel{Sidef}
S_i ~\doteq~ \bigl\{x \in \partial \Omega_i\,;~  \bfn(x) ~\hbox{is a singleton} \bigr\}.
\eeq
{Being $\Omega$ semiconvex and $y = x + r n(x)$ the only point in $\partial \Omega^r$ at distance $r$ from $x\in S_i$,} the open ball
$B_r\bigl(x+r\bfn(x)\bigr)$ cannot intersect $\Omega$.   As shown in Fig.~\ref{f:sc101},
left,
it now suffices to observe that 
$$\Omega_i~\cup~\bigcup_{x\in S_i}~B_r\bigl(x+r\bfn(x)\bigr)~\supseteq~ {(\clos\,\Omega_i)^\rho},$$
where $\rho$ is the radius at (\ref{domi}).
\endproof

{In particular, $\clos\,\Omega_i$ coincides with the connected components of $\Omega$ with positive measure.}

{From} the definition (\ref{Sidef}) it follows that  $\partial \Omega_i \setminus S_i$ is a set of  corner points, admitting multiple outer normals; hence it is countable. A further property
of points $x\in S_i$ is now described.

\begin{lemma}
\label{l:350}
For every point $x \in S_i$ there exists $\delta = \delta(x) > 0$ such that 
\bel{empty}
\bigl\{x + \rho  \bfn(x)\,;~r-\delta<\rho<r \bigr\}~ \subset ~{(\clos\,\Omega_i)}^r \setminus \clos (\Omega \setminus {\clos\,}\Omega_i)^r.
\eeq
\end{lemma}

\begin{figure}[ht]
\centerline{\hbox{\includegraphics[width=8cm]{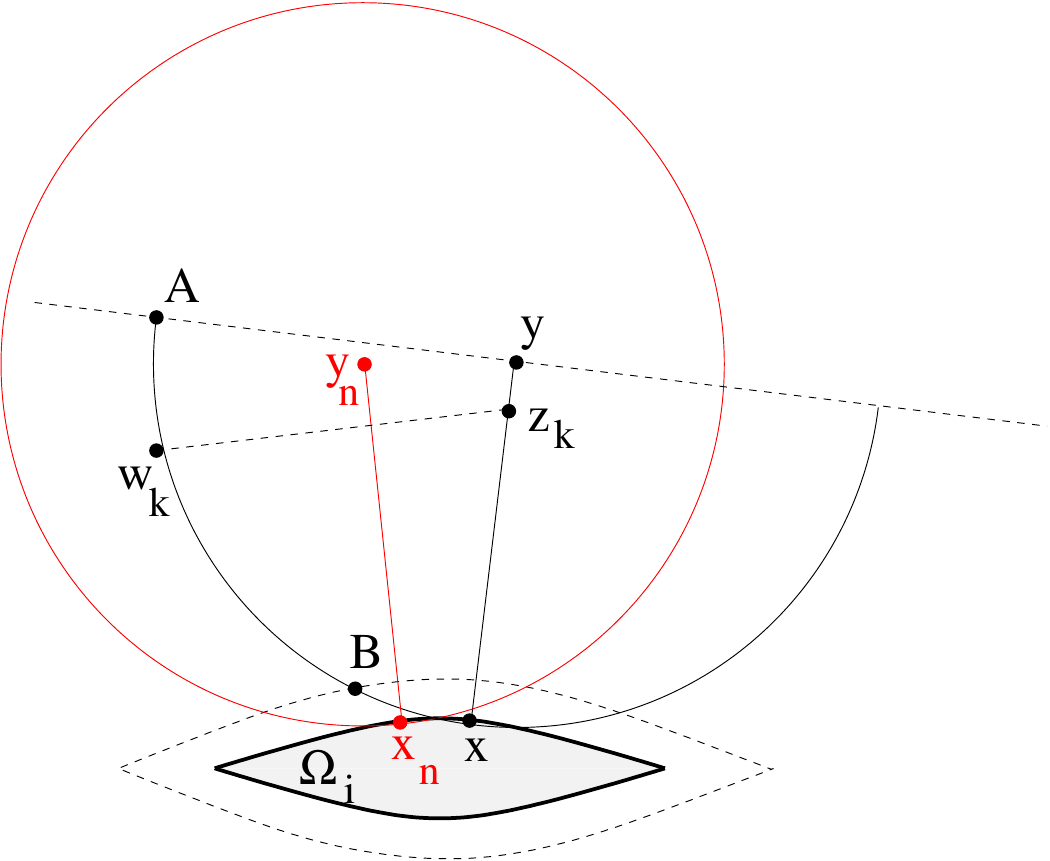}}}
\caption{\small Proving Lemma~\ref{l:350}. Here the arc of circumference $AB$
along the circumference centered at $y=x+r\bfn(x)$ with radius $r$ is entirely contained in the interior of the ball $B_r(y_n)$, for $n$ large enough. Therefore, for $k$ large, the point
$w_k$ cannot lie inside $\Omega$. }
\label{f:sc111}
\end{figure}

{\bf Proof.} With reference to Fig.~\ref{f:sc111}, consider a sequence of boundary points
$x_n\in S_i$ with $x_n\to x$.   Let $y= x+r \bfn(x)$ and $y_n=x_n+ r \bfn(x_n)$.
Assume that the conclusion of the Lemma fails.  Then there exists an increasing 
sequence $c_k\to r-$, and sequences of points $z_k, w_k$,  $k\geq 1$, such that 
$$z_k= x+ c_k\bfn(x),\qquad\qquad w_k\in \Omega\setminus {\clos\,\Omega_i}\,,\qquad |w_k-z_k|\leq r.$$
Since $z_k\to y$, by possibly taking a subsequence we conclude that 
$w_k\to \ov w\in \Omega$, with $|\ov w-y|=r$.   

By Lemma~\ref{l:340}, every point $w_k$ has uniformly positive distance from $\Omega_i$.  Hence the limit point $\ov w$ must lie on an arc $AB$ of the circumference
$\partial B_r(y)$ of length $<\pi r/2$.  However, this is impossible because
such arc is entirely contained in the open ball $B_r(y_n)$, for $n\geq 1$ large enough.
\endproof

\begin{lemma}
\label{l:34} Assume that {the interior $\inter\,\Omega$ of an optimal set $\Omega$} has infinitely many connected components.
Then there exists a sequence of components $\Omega_k$ such that
\begi
\item[(i)] ${\mathrm{diam}(\Omega_k)\to 0}${;}
\item[(ii)] {e}ach set $\Omega_k$ contains two corner points $x_k, x_k' {\in \partial \Omega_k}$, 
where the sets of outer normals
$\bfn(x_k)$, $\bfn(x_k')$ satisfy
\bel{nxn}
\bigl|\bfn(x_k)\bigr|~\to~\pi, \qquad\qquad \bigl|\bfn(x'_k)\bigr|~\to~\pi, \qquad \hbox{as}\quad {k} \to \infty{;}\eeq
\item {for $k \gg 1$, writing $\partial \Omega_k \setminus \{x_k,x_{k'}\}$ as the union of the two Lipschitz arcs $(\partial \Omega_k)^+,(\partial \Omega_k)^-$, only one of the two arcs may have nonempty intersection with $\partial \Omega_0$.}
\endi
\end{lemma}

\begin{figure}[ht]
\centerline{\hbox{\includegraphics[width=10cm]{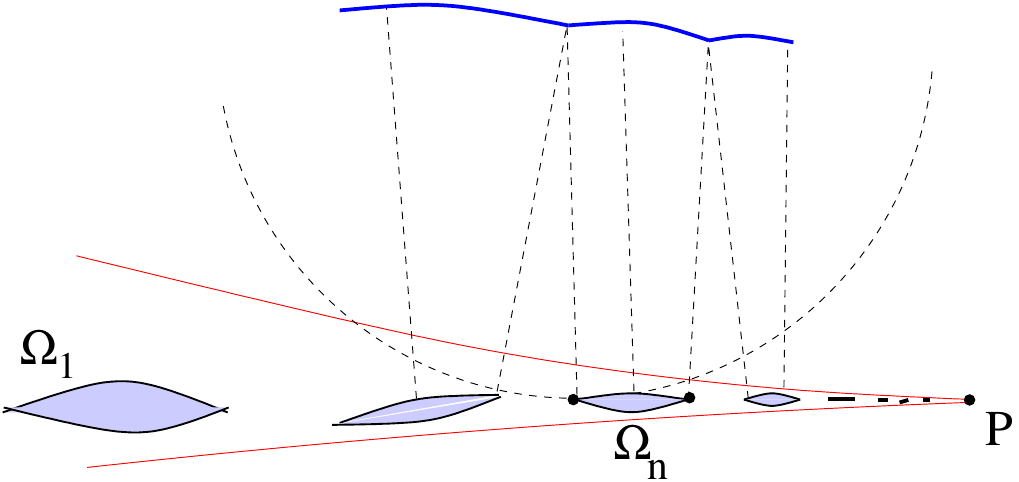}}}
\caption{\small The situation considered in Lemma~\ref{l:34}. }
\label{f:sc94}
\end{figure}

{\bf Proof.} {\bf 1.} 
As shown in Fig.~\ref{f:sc101}, left, consider any convex component  $\Omega_i\subset \inter\,\Omega$. 
Assume that, at each point $x\in \partial \Omega_i$, the set $\bfn(x)\subset~{\bbS^1}$ of outer normals covers an 
angle $\leq \theta$.   Then, by the duality relation $\Omega= (\Omega^r)^{-r}$, every other component
$\Omega_j$ must have distance from $\Omega {\setminus (\clos\,\Omega_i)}$ 
at least
$$\delta~\geq~|x'-x|~=~2 r\,\cos{\theta\over 2}\,.$$
\v
{\bf 2.} If $\inter\,\Omega$ contains infinitely many components, since each one of them has positive Lebesgue measure, there can be only countably many of them.   Since $\Omega$ is bounded, by taking a subsequence $(\Omega_k)_{k\geq 1}$ 
we can assume that their barycenters $b_k$ converge to some limit point $\bar x$.   By the previous step, 
there is a sequence of points $x_k\in \Omega_k$ where the sets of unit normal vectors have 1-dimensional measures
$\bigl|\bfn(x_k)\bigr|\to\pi$.  
As shown in Fig.~\ref{f:sc101}, right, call $\bfn_k = e^{i\theta_k}$ the central unit normal at $x_k$, so that 
the entire set of unit normals has the representation
$$\bfn(x_k) ~=~\left\{e^{i\theta}\,;~~\theta_k -\alpha_k \leq \theta\leq \theta_k +\alpha_k\right\},$$
for suitable angles $\alpha_k<\pi/2$.
By possibly taking a further subsequence, we can assume $\bfn_k\to \ov \bfn$. 
\v
{\bf 3.} Next, 
we claim that each set ${\partial} \Omega_k$ must also contain a second point $x'_k$, such that the 
sets of unit normal vectors $\bfn(x'_k)$ also satisfy $\bigl|\bfn(x_k')\bigr|\to\pi$.

Indeed, if the claim did not hold, we could find $\delta>0$ such that for each $k\geq 1$, the half circle
$$\Big\{y\in \R^2\,;~~|y-x_k|\leq \delta, ~~(y-x_k)\cdot \bfn(x_k)\leq 0\Big\}$$
does not intersect any other connected component besides $\Omega_k$.  This would exclude the existence of a limit 
point $\ov x$, providing a contradiction.
This establishes part (ii) of the {statement}.

Notice that, from the convergence $\bigl|\bfn(x_k)\bigr|\to\pi$, it follows
that the central unit normals $\bfn_k'$ at $x_n'$ satisfy $\bfn_k'\to - \ov \bfn$, as $k\to\infty$.

{If (iii) is false, then the limit point $\bar x$ belongs to $\partial \Omega_0$, and the suppporting cone $\R \bfn(\bar x)$ of $\Omega_0$ would have an opening of $\pi$, which is impossible by convexity if $\mathcal L^2(\Omega_0) > 0$.}
\v
{\bf 4.} 
To complete the proof, assume that (i) fails.   By possibly taking a subsequence, we can assume that 
${\mathrm{diam}}(\Omega_k)~=~|x_k-x'_k|~\geq ~\delta~>~0$ for every $k\geq 1$.   
Taking further subsequences, we obtain the convergence
$$x_k\to \ov x,\qquad x'_k\to \ov x ',\qquad\qquad |\ov x - \ov x'|\geq \delta.$$
This yields an obvious contradiction with the duality assumption: $\Omega= (\Omega^r)^{-r}$. 
\endproof

In connection with the optimization problem (\ref{argmin}),
the next lemma provides necessary conditions for a set $\Omega\subseteq {\Omega_0}$ to be optimal.
Namely,  every component $\Omega_i$ {of $\inter\,\Omega$} must have the same curvature radius $R$, at all 
points in the interior of ${\Omega_0}$, with the exception of countably many corner points described in Lemma~\ref{l:33}.
More precisely, define the  sets of boundary points
\bel{Sidef2}
{\tilde S_i} ~\doteq~ \big\{x \in \partial \Omega_i\,;~~ x + \bfn(x) r \subset \partial \Omega^r \big\}.\eeq
%
%
%
\begin{lemma}
\label{l:35}   Let {$\Omega$ be an optimal set for the problem (\ref{argmin}), and $\inter\,\Omega =\cup_i\Omega_i$}. Then the curvature of its boundary is constant in the set {$S = (\cup_i \tilde S_i) \setminus \partial {\Omega_0}$, and hence $\tilde S_i = S_i$, where $S_i$ is defined in (\ref{Sidef})}. 
\end{lemma}

{\bf Proof.} 
{The identity $S_i = \tilde S_i$ follows is all the points in $\tilde S_i$ have the same curvature, a conditions which implies that $\bfn(x)$ is singleton.}

This statement follows by removing a small area near a point $x'\in S$ where the boundary has a smaller curvature radius $R'$, and adding the same area near a point $x''\in S$ with larger curvature radius $R''$.
This can be done as long as $x''$ lies in the interior of ${\Omega_0}$.

More precisely, let $\Omega'$ be the set obtained from $\Omega$ by removing a small region near $x'$, as in  \eqref{E1h}, \eqref{E2h}.
In view of Lemma~\ref{l:350}, most of the points  removed from ${(\clos\,\Omega_i)}^r$
do not lie in the set $(\Omega \setminus {\clos\,\Omega}_i)^r$. 

By (\ref{Eper2}),  the change in the area of the $r$-neighborhood is estimated by
\begin{equation*}
{\mathcal L^2(\Omega^r \setminus (\Omega')^r)} ~=~ \bigg( 1 + \frac{r}{R'} + o(1) \bigg) {\mathcal L^2(\Omega \setminus \Omega')}.
\end{equation*}
Next,  let $\Omega''$ be the set obtained from $\Omega$ by adding a small region 
near $x''$, as in  (\ref{E'}).
By (\ref{E:p1convex}),  the change in the area of the $r$-neighborhood is estimated by
\begin{equation*}
{\mathcal L^2((\Omega'')^r \setminus \Omega^r)} ~=~ \bigg( 1 + \frac{r}{R''} + o(1) \bigg) {\mathcal L^2(\Omega'' \setminus \Omega)}.
\end{equation*}

By simultaneously performing the two modifications, we obtain a new set $\Tilde \Omega
\subset {\Omega_0}$, with ${\mathcal L^2(\Tilde \Omega)}={\mathcal L}^2 (\Omega)$ but $\mathcal L^2(\Tilde\Omega^r)<{\mathcal L^2(\Omega^r)}$,
against the optimality of $\Omega$.
\endproof

In the following, we shall denote by $R$ be the constant curvature radius, outside the corner points,  as in Lemma \ref{l:35}.

To prove that the configuration considered in Lemma~\ref{l:34} is not optimal, 
we study what happens if we remove one of the sets $\Omega_i$.   

As shown in Fig.~\ref{f:sc94}, let $x_i, x_i'\in {\partial \Omega_i}$ be the corner points where the set
of outer normals is large, say $\bigl|\bfn(x_i)\bigr|> \pi-\ve_0$, $\bigl|\bfn(x'_i)\bigr|> \pi-\ve_0$.
{By convexity, the boundary $\partial\Omega_i$ is the union of two parts, above and below the 
segment $[x_i, x_i']$, only one of which may have nonempty intersection with $\partial \Omega_0$ by Lemma \ref{l:34}}. To fix ideas, we assume that the upper boundary lies in the interior of ${\Omega_0}$.
The lower boundary may have nonzero intersection with $\partial {\Omega_0}$. 

\begin{figure}[ht]
\centerline{\hbox{\includegraphics[width=8cm]{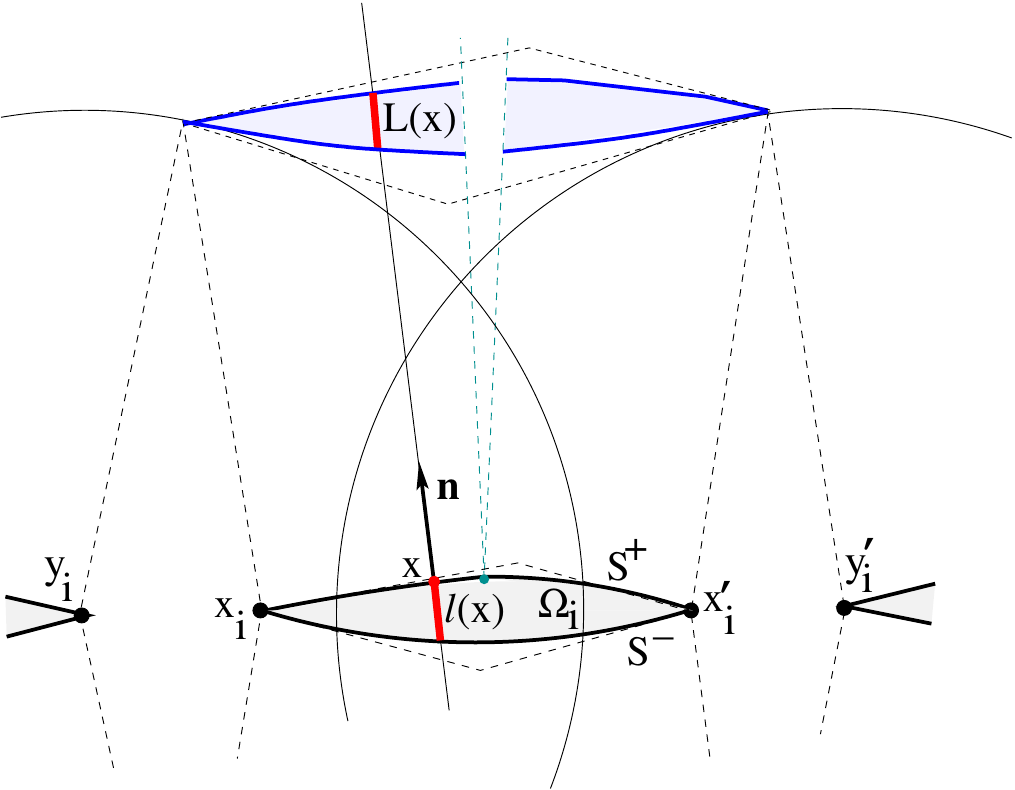}}}
\caption{\small The configuration considered at (\ref{Om7})-(\ref{Om8}).
Here $S^+_{i}\subset \partial\Omega_i$ is the upper part of the boundary of $\Omega_i$, excluding the points with multiple outer normals.
Removing the set ${\clos\,\Omega_i}$, all points in the upper shaded region are no longer in the set $(\Omega\setminus{\clos\,\Omega_i})^r$. These are points $y\in Y^+_{i}$ which project into $S^+_{i}$, but whose distance from $x_{i-1}'$ and from $x_{i+1}$ (and from all other components of $\Omega$) is $>r$.} 
\label{f:sc93}
\end{figure}

Call $S^+_{i} \subset \partial\Omega_i$ the set of points on the upper boundary 
having a unique outer normal, and consider the set of points projecting into $S^+_{i}$, namely
$$Y^+_{i}~\doteq~\Big\{y\in \R^2\,;~~d(y, {\Omega_i})=|y-x|~~\hbox{for some}~x\in S^+_{i}\Big\}.$$
We estimate how many of these points will not be in $(\Omega\setminus{\clos\,\Omega_i})^r$.
For this purpose, for each $x\in S^+_{i}$, consider the segments
\bel{lL1}\ell(x) \,\doteq\, \bigl\{x + c \bfn(x)\,;~c\in\R\bigr\} \cap \Omega_i\,,\qquad\qquad 
L (x) \,\doteq\, \ell(x) + {r}\, \bfn(x).
\eeq
We observe that, for every $x\in S^+_{i}$, $\bigl|L(x)\bigr|= \bigl|\ell(x)\bigr|$ and moreover
\bel{lL2}\ell(x)\,\subset\,\Omega_i\,,\qquad\qquad L(x)~\subset~\Omega^r\setminus (\Omega\setminus{\clos\,\Omega_i})^r.
\eeq
Indeed (see Fig.~\ref{f:sc93}),  let $x_{i-1}'\in {\partial}\Omega_{i-1}$ be the corner point
opposite to $x_i$ as in Lemma~\ref{l:34}. Similarly, let $x_{x+i}\in {\partial}\Omega_{i+1}$ be the corner point opposite to 
$x_i'$. 
We then have the identity 
$$Y^+_{i} \cap (\Omega\setminus {\clos\,}\Omega_i)^r~=~Y^+_{i} \cap \bigl(B_r(x_{i-1}') \cup
B_r(x_{i+1})\bigr).$$
Since $L(x)\cap \bigl(B_r(x_{i-1}') \cup
B_r(x_{i+1})\bigr)=\emptyset$, this yields the second relation in (\ref{lL2}).

By Lemma~\ref{l:35}, at every point $x\in S_{i}$, the curvature of $\partial \Omega_i$ is constantly equal to $R$.

Being the curve $\partial \Omega_i$ convex, the map
\begin{equation*}
S_i^+ \times \R \ni x,c \mapsto T(x,c) = x + c \n(x) \in \R^2
\end{equation*}
is BV on the rectifiable set $S_i^+ \times \R$, and its a.c. Jacobian is
\begin{equation*}
J(x,c) = \bigg| 1 + \frac{c}{R} \bigg|.
\end{equation*}
This allows us to compute by the area formula
\begin{equation}
\label{Om7}
\begin{split}
{\mathcal L^2(\Omega_i)}~&\leq~\int_{\Omega_i} \mathcal H^0(T^{-1}(y)) \mathcal L^2(dy) = \int_{S^+_i \times \R^- \cap T^{-1}(\Omega_i)} J(x,c) \mathcal H^1(dx) \mathcal L^1(dc) \leq \int_{S_i^+} |\ell(x)| \mathcal H^1(dx),
\end{split}
\end{equation}
where we observe that, being the set convex, $\Omega_i \subset T(S_i^+ \times \R^-)$. Similarly
\bel{Om8} \bega{rl}
{\mathcal L^2 \big( (\Omega\setminus\clos\,\Omega_i)^r \cap Y^+_i \big)}&\ds=~\mathcal H^2 \big( S^+_i \times \R^+ \cap T^{-1}((\Omega\setminus\clos\,\Omega_i)^r \cap Y^+_i) \big) \\[4mm]
&\ds \geq~\int_{S^+_i} \bigg[ \int_{r - |L(x)|}^{r} \bigg( 1 + \frac{z}{R} \bigg) \, dz \bigg] \mathcal H^1(dx) \\[4mm] 
&\ds =~\int_{S^+_i} \bigg( r + \frac{r^2}{R} - (r-L(x)) - \frac{(r-L(x))^2}{R} \bigg) \, \mathcal H^1(dx) \\[4mm] 
&\ds = ~\int_{S^+_i} \bigg( {R+r\over R} - {\bigl|\ell(x)\bigr|\over 2R} \bigg) \bigl| \ell(x)\bigr|\, dx.
\enda\eeq 
Combining the above inequalities, we conclude
\bel{Om9} {
{\mathcal L^2 \bigl((\Omega\setminus{\clos\,\Omega_i})^r\cap Y^+_{i}\bigr)}\over
{\mathcal L^2(\Omega_i)}}~\geq~1+{r\over R } - o(1),\eeq
where $o(1)$ is a quantity that approaches zero as $\sup_x\bigl| \ell(x)\bigr|\to 0$.

A similar estimate can be performed for the set of points $y\in (\Omega\setminus {\clos\,\Omega_i})^r\cap Y^-_{i}$  which project onto the lower boundary $S^-_{i} \subset\partial \Omega_i$.  However, if this lower boundary touches $\partial {\Omega_0}$,
we have no lower bound on its curvature. With the same computations as above can only obtain the weaker estimate
\bel{Om11} {
{\mathcal L^2((\Omega\setminus{\Omega_i})^r\cap Y^-_{i})}\over
{\mathcal L^2(\Omega_i)}}~\geq~1.\eeq
This enough to conclude
the non-optimality of $\Omega$.  Indeed, when $\Omega_i$ is small enough, (\ref{Om9}) and (\ref{Om11}) 
together yield
$${{\mathcal L^2((\Omega\setminus{\clos\,\Omega_i})^r}\over
{\mathcal L^2(\Omega_i)}}~\geq~{3\over 2} + {r\over R}\,,$$
providing a contradiction.
We thus have
\v

\begin{lemma}
\label{l:36}
The optimal set $\Omega$ can have at most finitely many connected components.
\end{lemma}

Indeed, if there are infinitely many components, then there exists an accumulation point and thus we are in the situation of Lemma \ref{l:34} for arbitrary small sets: deleting one of these very small components
and transferring
its mass along the boundary of another set $\Omega_j$ with curvature $R$, the total area ${\mathcal L^2(\Omega^r)}$ will decrease.

{In particular, the connected components of $\Omega$ are the closure of the connected components $\Omega_i$, $i=1,\dots,N$: with a slight abuse of notation, we will use the notation $\Omega_i$ for the components of $\Omega$.}

\subsection{Completion of the proof of Theorem~\ref{t:31}.}
It now remains to prove that a set $\Omega= \cup_{i=1}^N\Omega_i$, with a finite number $N\geq 2$ of connected components, 
is not optimal. 
W.l.o.g., we can assume that ${\Omega_0} = \conv \,\Omega$. We will show that it is possible to 
rigidly move each component, so that
\bel{omit}
\Omega_i (t) ~=~t(z-z_i)+ \Omega_i,
\eeq
in such a way that the area of the $r$-neighborhood
$${\mathcal L^2(\Omega(t)^r)}~=~\sum_i {\mathcal L^2(\Omega_i(t)^r)} - \sum_{i \not= j} {\mathcal L^2 \big( \Omega_i(t)^r \cap \Omega_j(t)^r \big)}$$
is strictly decreasing.  Here the points $z\in {\Omega_0}$ and $z_i\in \Omega_i$  must be carefully chosen, so that 
all the sets $\Omega_i(t)$
remain inside ${\Omega_0}$, for $t\in [0,\ve]$ with $\ve>0$ small enough. {In the above formula we have used \eqref{Equa:triple_intersection}.}

If $\Omega= \cup_i\Omega_i$ is an optimal set,  
the time derivative of the  area ${\mathcal L^2(\Omega(t)^r)}$ is computed as follows.
Call ${\cal A}$ the set of all couples $(i,j)$ such that $i\not= j$ and 
$$d_{ij}~\doteq ~\dist(\Omega_i,\Omega_j)~=~|x_i-x_j|~<~2r.$$
Here $x_i, x_j$ are the points considered in Lemma \ref{l:34}.
At time $t=0$, assuming that $z_i\in \Omega_i$ for all $i=1,\ldots,N$, the time derivative is computed by
\bel{dar}
\frac{d}{dt} \mathcal L^2 \bigg( \bigcup_i \Omega_i(t)^r \bigg)_{t=0} ~= ~ \sum_{(i,j)\in \A} 2 \sqrt{r^2 - \frac{d_{ij}^2}{4}} ~(x_i - x_j) \cdot (z_j-z_i)~\leq~0.
\eeq
We can always assume that $\Omega^r=\cup_i \Omega_i^r$ is connected, otherwise the non-optimality is trivial{: indeed, with the same ideas here below, otherwise there is a component of $\Omega^r$ which can be moved freely inside $\Omega_0$ until it superimpose to another component of $\Omega^r$}.  Notice that this assumption 
implies that, for every $i$, there exists $j\not=i$ such that $\Omega_{i}^r \cap \Omega_{j}^r \not= \emptyset$. As a consequence, we can assume that the inequality (\ref{dar}) 
is strict.

\begin{figure}[ht]
\centerline{\hbox{\includegraphics[width=8cm]{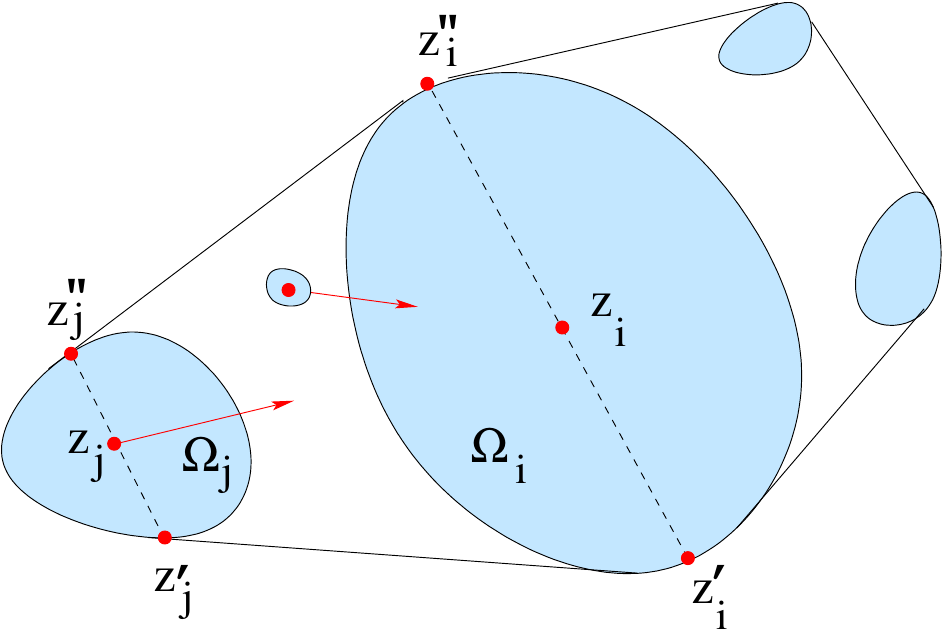}}
}
\caption{\small  Proving that a set $\Omega= \bigcup_i\Omega_i \subset {\Omega_0}$ with 
finitely many connected components cannot be optimal.    The case considered in Lemma~\ref{l:37}.}
\label{f:sc105}
\end{figure}

If $\Omega$ has at least two components,
the first lemma rules out the existence of an $\Omega_i$ which cannot be translated inside $ {\Omega_0} $ (see Fig.~\ref{f:sc105}).

\begin{lemma}
\label{l:37} Let $\Omega\subset {\Omega_0}$ be an optimal set.
If $\Omega$ has more than one component, then for every component $\Omega_i$ the set $\partial \Omega_i \cap \partial  {\Omega_0} $ is a connected arc whose normal vectors span an angle $< \pi$.
\end{lemma}

{\bf Proof.} {\bf 1.} For a component $\Omega_i$ such that  $\partial \Omega_i \cap \partial  {\Omega_0}=\emptyset $, the spanned angle is zero and the conclusion is trivial.
\v
{\bf 2.}
Next,
assume that there is a component  $\Omega_i$ such that $\partial  {\Omega_0}  \setminus \partial \Omega_i$ contains an arc $\arc{z_i'z''_i}$ such that the angle spanned by $\bfn(x)$, $x \in \arc{z_i'z_i''}$, is $\leq \pi$. We consider only the components  $\Omega_j$ contained in the part of $ {\Omega_0}  \setminus \Omega_i$ whose boundary contains the arc $\arc{z_i'z_i''}$.  There must be at least one such set, because $ {\Omega_0}  = \conv\, \Omega$. 
To obtain a contradiction, we shall move only these sets.

If $\partial \Omega_j \cap \partial  {\Omega_0}  = {\arc{z_j'z_j''}}$, we consider the point  $z_j = \frac{z_j' + z_j''}{2} \in \Omega_j$.  
Otherwise, if $\partial \Omega_j \cap \partial  {\Omega_0} = \emptyset$, we take the barycenter:
$z_j = \avint_{\strut \Omega_j} x\, dx$. 
In addition, we set  $z = z_i = \frac{z_i' + z_i''}{2}$.
Defining $\Omega_i(t)= \Omega_i$ while 
$\Omega_j(t) = \Omega_j + t(z-z_j)$, we check that $\Omega_j(t) \subset  {\Omega_0} $  for $t>0 $ small.  Moreover, (\ref{dar}) is satisfied as a strict inequality.
Hence $\Omega$ is not optimal.
\endproof
\v
If $\Omega$ contains  more than one component, 
by Lemma~\ref{l:37}, for every component $\Omega_i$ the intersection $\partial \Omega_i \cap \partial  {\Omega_0} $ either is empty, or is a connected arc spanning an angle $< \pi$.
The next lemma rules out this remaining possibility.

\begin{figure}[ht]
\centerline{\hbox{\includegraphics[width=11cm]{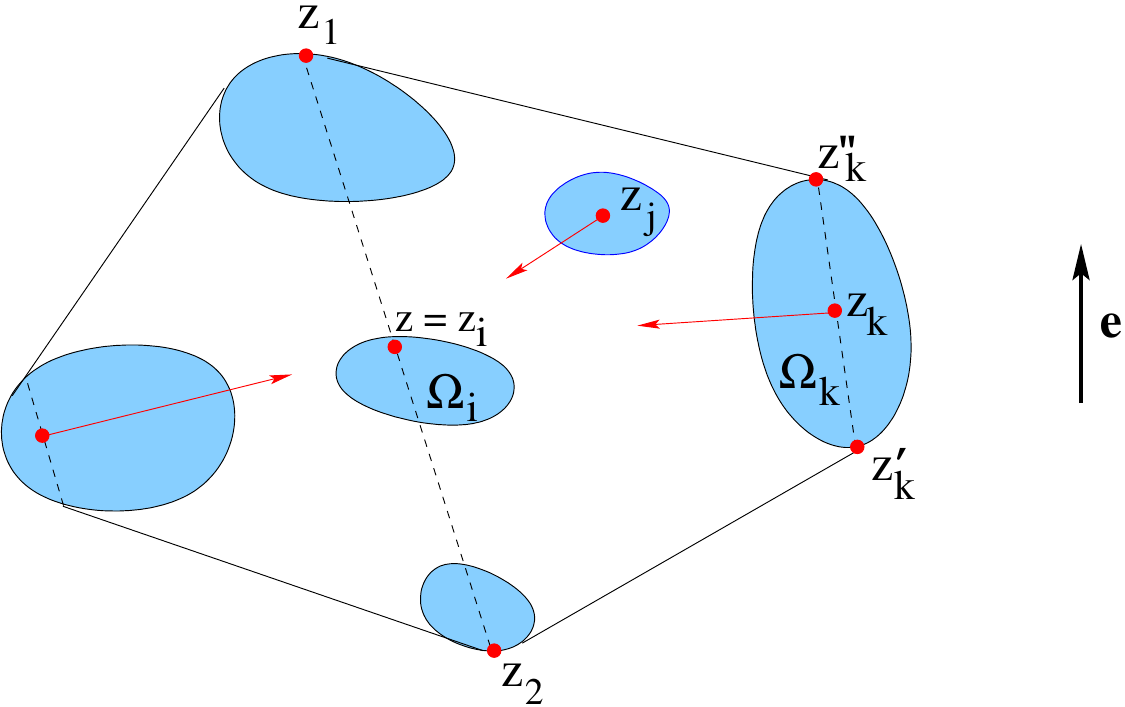}}
}
\caption{\small The configuration  considered in Lemma~\ref{l:38}.}
\label{f:sc109}
\end{figure}

\begin{lemma}
\label{l:38}
If $\Omega\subset {\Omega_0}$ has more that one component, and every component intersects $\partial  {\Omega_0} $ in an arc of opening $<\pi$ (or does not touch $\partial  {\Omega_0} $ at all),  then $\Omega$ is not optimal.
\end{lemma}


{\bf Proof.} For every component $\Omega_i$, by assumption the intersection
$\partial \Omega_i \cap \partial  {\Omega_0} $ is a connected arc ${ \arc{z_i'z_i''}}$ of opening $< \pi$, or it is empty (see Fig.~\ref{f:sc109}). 

Fix a unit vector $\bfe\in \bbS^1$, and choose two points $z_1,z_2 \in \partial  {\Omega_0} $ such that $\bfe \in \bfn(z_1)$, $\bfe\in -\bfn(z_2)$. Define $z = \frac{z_1+z_2}{2}$.  Then, for every 
component $\Omega_i$,  if $z \in \Omega_i$   we define  $z_i \doteq  z$.
On the other hand, if   $z \notin \Omega_i$, we consider two cases.
\begi
\item If $\partial \Omega_i \cap \partial  {\Omega_0}  ={ \arc{z_i'z_i''}}$, we take the mid-point:  $z_i \doteq \frac{z_i' + z_i''}{2}$.
\item If $\partial \Omega_i \cap \partial  {\Omega_0}  = \emptyset$, we take the barycenter: $z_i = \avint_{\strut\Omega_i} x \,dx$. 
\endi
By elementary geometry, we check that, for $t>0$ small, all the sets $\Omega_i(t)$ defined at \eqref{omit} remain inside $ {\Omega_0} $.  Again, the relation (\ref{dar}) is satisfied as a strict
inequality.
Hence $\Omega$ is not optimal.\endproof

Combining the two above lemmas, we achieve  the proof of Theorem~\ref{t:31}.

\v

\section{Further properties of the optimal set}
\label{s:4}
\setcounter{equation}{0}
Having proved that the optimal set for the problem (\ref{op1}) is convex, by the area formula
$${\cal L}^2(\Omega^r) ~=~{\cal L}^2(\Omega)+ r\,{\cal H}^1(\partial \Omega)+ \pi r^2$$
it is clear 
that the same set $\Omega$  is {an} optimal set for the constrained
isoperimetric problem (\ref{op2}).   This second problem has been studied in
\cite{SZ}.  Recalling the definition of the sets $\Hat\Omega$ 
 at (\ref{HOdef}), 
we collect here the main results:
\begin{theorem}
\label{t:41} Let $\Omega_0\subset\R^2$ be a compact, convex set. Then, for every
$0<a\leq {\mathcal L^2(\Omega_0)}$ the constrained  minimization problem (\ref{op2}) has a solution. 
Moreover, the  following holds.
\begi

\item[(i)] The optimal set  $\Tilde\Omega=\Tilde \Omega(\Omega_0 ,a)$ is convex. Moreover  its boundary $\partial \Tilde\Omega$ has curvature $\kappa$ which is constant and maximal along each connected arc in $\partial \Tilde \Omega \setminus \partial \Omega_0 $.

\item[(ii)] Conversely, any convex set $\Tilde\Omega\subseteq \Omega_0$, with 
${\mathcal L^2(\Tilde \Omega)}=a$ and such that the curvature is maximal and constant 
along $\partial \Tilde \Omega \setminus \partial \Omega_0 $,  is an optimal solution to (\ref{op1}).

\item[(iii)] When $0<a\leq\pi \ov  R^2$, with $\ov R$ being the inner radius in (\ref{ovrdef}),
the optimal solutions are precisely the  balls 
$B_\rho(x) \subseteq \Omega_0 $ with radius $\rho=\sqrt{a/\pi}$.

\item[(iv)] When $\pi \ov  R^2 < a \leq {\mathcal L^2(\Hat \Omega(\Omega_0 ,\ov  R))}$, any 
optimal solution is contained in $\Hat \Omega(\Omega_0 ,\ov  R)$, and coincides with the convex closure of two balls of radius $\ov R$. 

\item[(v)] For ${\mathcal L^2(\Hat \Omega(\Omega_0 ,\ov  R))}\leq a \leq \mathcal L^2(\Omega_0)$, 
the optimal solution is unique. Indeed,
there exists a unique $\rho\in \,]0, \ov R]$ such that
\begin{equation*}
\Tilde \Omega(\Omega_0 ,a) ~=~ \Hat \Omega(\Omega_0 ,\rho).
\end{equation*}

\item[(vi)] For  $0<a<\mathcal L^2(\Omega_0)$, one has
\bel{bcurv}
\frac{d}{da} {\mathcal H^1(\partial \Tilde \Omega(\Omega_0 ,a))}~ =~ \kappa,\eeq
where $\kappa$ is the maximal curvature of the boundary $\partial \Tilde\Omega$.

\item[(vii)] 
The map $a \mapsto {\mathcal H^1(\partial \Tilde \Omega(\Omega_0 ,a))}$ is monotone increasing.
\endi
\end{theorem}

{\bf Proof.} {\bf 1.}
The existence and the convexity of solutions were the main results proved in Theorem~3.31
of \cite{SZ}, together with the properties stated in (i) and (vii).  
Parts (iii),  (iv), (v) are proved in Theorem~3.32 of \cite{SZ}.
\v  
{\bf 2.}
To prove (ii), let $\Omega \subset \Omega_0$ be such that ${\mathcal L}^2 (\Omega) = a$ and the curvature $\kappa$ is constant and maximal in the set $\partial \Omega \setminus \partial \Omega_0$.   
Any relatively open connected component  of $\partial \Omega \setminus \partial \Omega_0$ is thus an arc of a circle of radius $\rho = 1/\kappa > 0$.  

{Being the of every point of the boundary $\kappa(x) \leq 1/\rho$}, it follows that $\Omega$ has the internal ball property. Namely, for every $x \in \Omega$ there exists a ball $B_\rho(y)$ such that $x \in {\hbox{clos}}\, B_\rho(y) \subset \Omega$. Therefore, $\Omega$ is the closure of a union of balls  of radius $\rho$. 

We use the following observation.  Let $\Omega$ be a compact convex set set whose boundary 
has curvature $\kappa(x) \leq 1/\rho$ at every point $x\in\partial \Omega$. 
If the boundary is tangent to a ball $B_\rho(y)$ at two points $z_1,z_2 \in \partial B_\rho(y)$, then  $\Omega$ must contain the arc ${\arc{z_1z_2}}\subset \partial B_\rho(y)$ with minimal length (in case where both arcs have equal length $\pi \rho$, it must contain at least one of the two arcs). \\
If $\Omega$ contains an arc with length  strictly greater than $\pi \rho$ (i.e., spanning an angle $>\pi$), then it follows that $\Omega$ itself is a ball of radius $\rho$ {(case (iii) of the statement)}. 
\\
If an arc ${\arc{z_1z_2}}$ has length  exactly $\pi \rho$ (i.e., it spans an angle $=\pi$), then the set $\Omega$ must be the convex hull of two balls of radius $\rho$. Namely: $\Omega=\clos \Big(\bigcup_{\alpha \in [0,1]} B_\rho\bigl((1-\alpha)x_1 + \alpha x_2\bigr)\Big)$ for some $x_1,x_2$. Here we can take $x_1$ to be the center of the ball whose boundary contains the arc ${\arc{z_1z_2}}$,  while $x_2$ is the center of the inner ball tangent to $\partial\Omega$ at the  furthest point from $x_1$. Hence $\rho = \ov R$. {This corresponds to case (iv) of the statement.}

If every arc of radius $\rho$ contained in $\partial\Omega$  has length $< \pi \rho$ (i.e., it
spans an angle  $< \pi$), we claim that
$\Omega = \Hat \Omega(\Omega_0,\rho)$ {(case (v) of the statement)}.
In other words, it is impossible to enlarge $\Omega$ by adding other balls $B_\rho(y)$ of radius $\rho$ contained in $\Omega_0$. To prove the claim, two cases are considered:
\begin{itemize}
    \item If the ball $B_\rho(y)$ we add does  not intersect $\Omega$, then by convexity one arc must be $\geq \pi$ (which contradicts the assumptions here), otherwise there are two points with distance $>2\rho$ connected by an arc of radius $\rho$, which is impossible.
    \item On the other hand, if $B_\rho(y)\cap \partial\Omega\not=\emptyset$, then  the intersection must be an arc ${\arc{z_1z_2}}$ of curvature $\rho$.  But since the opening of this arc is $<\pi$,  {$B_\rho(y)$} must be already in $\Omega$; otherwise $B_\rho(y)$ is not {a subset of} $\Omega_0$. 
\end{itemize}

\begin{figure}[ht]
\centerline{\hbox{\includegraphics[width=7cm]{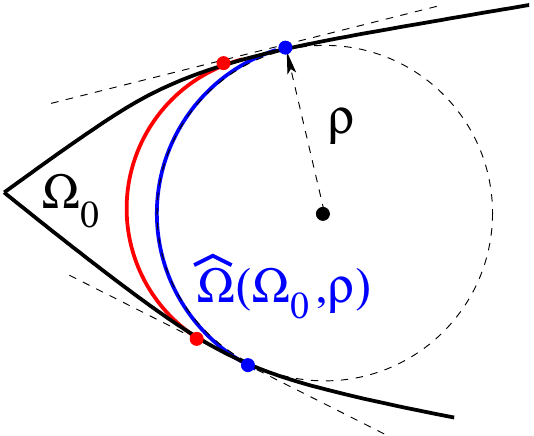}}
}
\caption{\small  Computing the increase in the area of  in $\Hat \Omega(\Omega_0, \rho)$, as the 
radius $\rho$ is reduced.}
\label{f:sc110}
\end{figure}

{\bf 3.} To prove (vi), let $\theta_i$ be the angle covered by the free arc $L_i\subset \partial\Omega\cap \Omega_0$ of maximal curvature, $i \geq 1$. Since by definition
\begin{equation*}
    \partial \Omega \setminus \partial \Omega_0 ~= ~\bigcup_i L_i, \qquad {\theta_i} = \frac{{{\mathcal H}^1}(L_i)}{\rho},
\end{equation*}

it is enough to study the variation of area and perimeter about a single arc $L_i$, spanning an angle ${\theta_i}$, which has constant curvature $\kappa = 1/\rho$. 
Moreover, the only non elementary case is when the opening of the arc is $< \pi$.  
In this case we will use the dependence w.r.t. the maximal radius of curvature $\rho$ of \eqref{HOdef} for $0<\rho <\ov R$.\\
Up to a rigid change of coordinates, about every arc $L_i$ we are in the situation shown in Fig.~\ref{f:sc110}. {If} ${\mathcal H^1(L_i(\rho))}$ is the length of the arc ${\arc{z_1z_2}}\subset \partial \Hat \Omega(\Omega_0,\rho)$, one gets
\bel{42}
    \bega{rl}
        {\mathcal H^1(L_i(\rho))}~ <~ {\mathcal H^1(L_i(\rho-h))} &=~ (\rho-h) \theta_i + 2 h \tan \frac{\theta_i}{2} + o(h) \\[2mm]
        &=~ {\mathcal H^1(L_i(\rho))} + h \bigg( 2 \tan \frac{\theta_i}{2}  - \theta_i \bigg) + o(h) \\[2mm]
        &\leq~ {\mathcal H^1(L_i(\rho))} + h \bigg( 2 \tan\frac{\theta_i}{2} - \theta_i \bigg),
    \enda
\eeq
where $\theta_i = {\mathcal H^1(L_i)}/\rho$.  Here we observed that $o(h) \leq 0$ by convexity.  
Note that the case of maximal growth occurs when $\partial \Hat \Omega(\Omega_0,\rho)$ coincides with the tangent cone. In particular, the map $\rho \mapsto {\mathcal H^1(L_i(\rho))}$ is right differentiable with derivative bounded by
\begin{equation*}
    2 \tan \bigg( \frac{\theta_i}{2} \bigg) - \theta_i.
\end{equation*}
The above arguments yield the estimate
\begin{equation*}
    {\mathcal H^1(\partial \Hat \Omega(\Omega_0,\rho))} ~< ~{\mathcal H^1(\partial \Hat \Omega(\Omega_0,\rho-h))}~ \leq ~{\mathcal H^1(\partial \Hat \Omega(\Omega_0,\rho))}
     + h \sum_i \bigg( 2 \tan \frac{\theta_i}{2}  - \theta_i \bigg).
\end{equation*}
Therefore the function $\rho \mapsto {\mathcal H^1(\partial \Hat \Omega(\Omega_0,\rho))}$ is Lipschitz for $\rho < \ov R$ and strictly decreasing w.r.t. $\rho$.  For every $k\geq 1$, taking the derivative 
for a.e.~$0 <\rho<\ov R$ we obtain
\begin{equation*}
    \begin{split}
		\sum_{i \leq k} \bigg( 2 \tan \frac{\theta_i}{2}  - \theta_i \bigg)~ &= {\sum_{i \leq k} - \frac{d}{d\rho} {\mathcal H^1(L_i(\rho))} } \\[4mm]
&\leq~ - \frac{d}{d\rho} {\mathcal H^1(\partial \Hat \Omega(\Omega_0,\rho))} ~\leq ~\sum_i \bigg( 2 \tan  \frac{\theta_i}{2} - \theta_i \bigg).
	\end{split}
\end{equation*}
Since the series is convergent,  for a.e.~$\rho$ we conclude that
\begin{equation*}
 - \frac{d}{d\rho} {\mathcal H^1(\partial \Hat \Omega(\Omega_0,\rho))} ~= ~\sum_i \left( 2 \tan \frac{\theta_i}{2} - \theta_i \right).
\end{equation*}

The same computation can be done for the area: the variation of area $A_i(\rho)$ inside each region of opening $\theta_i$ is computed by
\bel{atan}
\bega{rl}
     {\mathcal L^2(A_i(\rho))}~ <~ {\mathcal L^2(A_i(\rho-h))} &=~\ds {\mathcal L^2(A_i(\rho))} + \bigl( \rho^2 - (\rho - h)^2 \bigr) \bigg( \tan  \frac{\theta_i}{2} - \frac{\theta_i}{2} \bigg) + o(h)
     \enda
\eeq
with $o(h)\leq 0$. 
Arguing as before, for a.e.~$0<\rho<\ov R$ we obtain
\begin{equation*}
    - \frac{d}{d\rho} {\mathcal L^2(\Hat \Omega(\Omega_0,\rho))} ~=~ \rho \sum_i 
    \bigg( {2} \tan \frac{\theta_i}{2}  - \theta_i \bigg) ~= ~ - \rho \frac{d}{d\rho} {\mathcal H^1(\partial \Hat \Omega(\Omega_0,\rho))}.
\end{equation*}
This yields  \eqref{bcurv}, for a.e.~$\rho$. 
The final observation is that $\rho \mapsto {\mathcal L^2(\Hat \Omega(\Omega_0,\rho))}$ is strictly decreasing and continuous, as seen by \eqref{atan}. Therefore the inverse map
\begin{equation*}
    {\mathcal L^2(\Hat \Omega(\Omega_0,\rho))} \mapsto \rho
\end{equation*}
is continuous, thus implying that \eqref{bcurv} holds for all $\rho \in \,] 0,\ov R[\,$.
\endproof

\begin{figure}[ht]
\centerline{\hbox{\includegraphics[width=12cm]{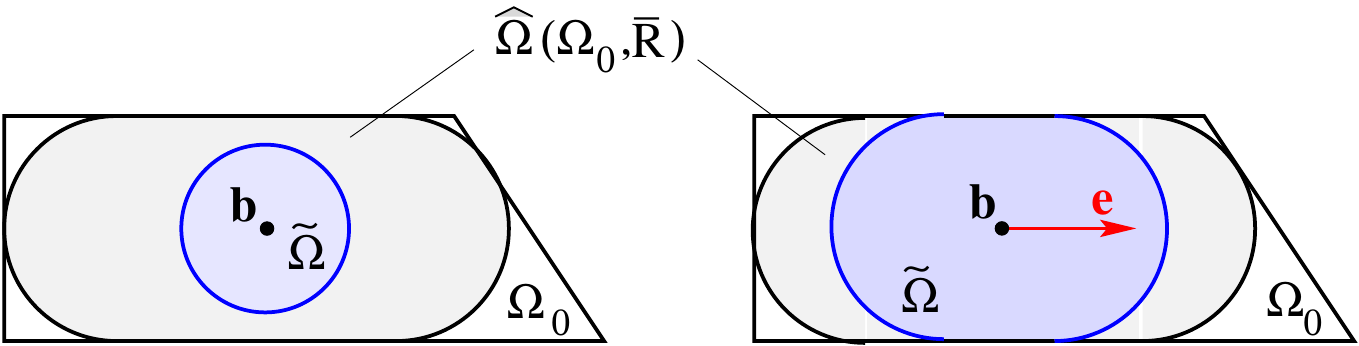}}
}
\caption{\small  Optimal sets $\Tilde \Omega = \Tilde\Omega(\Omega_0,a)$
in the cases considered at  (iii) and at (iv) of Theorem~\ref{t:31}, respectively.   Here
$\Omega_0$ is a trapezoid. Note that in both of these cases the optimal sets are not unique.
Uniqueness can be achieved by imposing that the barycenter $\bfb$ of $\Tilde\Omega$ coincide
with the barycenter of $\Hat \Omega(\Omega_0, \ov R)$.
}
\label{f:sc80}
\end{figure}

\begin{remark}\label{r:31} {\rm
We observe that the minimizer $\Tilde\Omega=\Tilde\Omega(\Omega_0 ,a)$ is not uniquely determined  when
\bel{asmall}
a ~<~ {\mathcal L^2(\Hat \Omega(\Omega_0,\ov R))}.
\eeq
As shown in Fig.~\ref{f:sc80}, one
can remove this ambiguity and single out a unique set $\Tilde\Omega$ by considering the barycenter
$\bfb$ of  $\Hat \Omega(\Omega_0 ,\ov  R)$, and imposing that the barycenter of $\Tilde\Omega$ also coincide with $\bfb$.  

More precisely, 
for a suitable unit vector $\bfe$ and $ \ell^*\geq 0$, we  have
the representation
$$\Hat \Omega(\Omega_0 ,\ov R) ~=~  \bigcup_{|h|\leq \ell^*}  B_{\ov R}(\bfb + h \bfe).
$$
We can then  define
\bel{TOM}
\Tilde\Omega(\Omega_0 ,a) ~\doteq~ \left\{\bega{cl}
B_{\sqrt{a/\pi}}(\bfb) \qquad & \hbox{if}\quad 0 \leq a \leq \pi \ov R^2,\\[4mm]
\ds\bigcup_{|h|\leq \ell} B_{\ov R}(\bfb+h \bfe)  \qquad & \begin{split} &\hbox{if}\quad \pi \ov R^2 <a <{\mathcal L^2(\Hat \Omega(\Omega_0 ,\ov R))}, \\[-1mm]
&\hbox{choosing $\ell$ so that} ~2 \ov  R \ell + \pi \ov  R^2 =a. \end{split}
\enda\right.
\eeq
 On the other hand, when ${\mathcal L^2(\Hat \Omega(\Omega_0,\ov R))} 
\,\leq\,a\,\leq\,{\mathcal L^2(\Omega_0)}$, the optimal set $\Tilde\Omega(\Omega_0 ,a) $ is already uniquely 
determined.}
\end{remark}

\begin{lemma}
\label{l:39}
With the definition (\ref{TOM}), the following holds:
\begi
\item[(1)] ~$a_1 < a_2$ implies $\Tilde \Omega(\Omega_0 ,a_1) \subset\Tilde \Omega(\Omega_0 ,a_2)$.
\item[(2)] ~$\Tilde  \Omega^r(\Omega_0,a)~ =~ \Tilde \Omega\Big(\Omega_0 ^r,\, {\mathcal L^2(\Tilde \Omega^r(\Omega_0 ,a))}\Big)$.
\endi
\end{lemma}

{\bf Proof.}
Part (1) of the lemma is obvious.  To prove part (2), we first observe that,
if $a \leq {\mathcal L^2(\Hat \Omega(\Omega_0 ,\ov  R))}$, 
the  result is a simple consequence of the identity
\begin{equation*}
\Hat \Omega(\Omega_0 ^r,\ov  R+r) ~=~ \Hat \Omega^r(\Omega_0 ,\ov  R).
\end{equation*}

To handle the remaining case where
$a > {\mathcal L^2(\Hat \Omega(\Omega_0 ,\ov  R))}$, we will prove that the set 
$\Tilde\Omega^r(\Omega_0, a)$, i.e., the neighborhood of radius $r$ around
$\Tilde\Omega(\Omega_0, a)$, provides an optimal solution to the problem
\begin{equation}
\label{op1r} \hbox{minimize:}~~~ {\mathcal L^2(\Omega^r)},\qquad\hbox{subject to}\quad\Omega 
\subseteq \Omega_0^r ,
\quad  {\mathcal L}^2 (\Omega) = {\mathcal L^2(\Tilde\Omega^r(\Omega_0, a))}.
\end{equation}

Toward this goal,
we claim that $\Tilde\Omega^r(\Omega_0, a)$ satisfies the
conditions stated in part (ii) of  Lemma~\ref{l:31}.
Indeed,
 consider any point $x\in \partial \Tilde\Omega(\Omega_0, a)$, and call 
$\bfn(x)$ the unit outer normal.
Then the curvature of $\Tilde \Omega^r(\Omega,a)$ at 
the boundary point $x + r \mathbf n(x) \in \partial \Tilde \Omega^r(\Omega_0 ,a)$ is computed by
\bel{curv}
\kappa\Big(\Tilde \Omega^r(\Omega_0,a),\, x + r \mathbf n(x)\Big) ~=~ \frac{\kappa\bigl(\Tilde \Omega(\Omega_0,a),x\bigr)}{1 + r \kappa\bigl(\Tilde \Omega(\Omega_0,a),x\bigr)}.
\eeq
Since by assumption the curvature $\kappa\bigl(\Tilde \Omega(\Omega_0,a),x\bigr)$
is constant and maximal at points $x\in \partial \Tilde \Omega(\Omega_0,a)\setminus \Omega_0$,
by (\ref{curv}) the curvature $\kappa\Big(\Tilde \Omega^r(\Omega_0,a),\, x + r \mathbf n(x)\Big)$ is constant and maximal at points  $ x + r \bfn(x)\in \partial \Tilde \Omega^r(\Omega_0,a)\setminus \Omega_0^r$.
Using (ii) in  Lemma~\ref{l:31},   we thus obtain the optimality of the set $\Tilde \Omega^r(\Omega_0,a)$.  

Finally, we notice the implication
\begin{equation*}
a \,>\, {\mathcal L^2(\Hat \Omega(\Omega_0 ,\ov  R))} \qquad \Longrightarrow \qquad {\mathcal L^2(\Tilde \Omega^r(\Omega_0 ,a))} \,>\,
 {\mathcal L^2(\Hat \Omega(\Omega_0 ^r,\ov  R+r))}.
\end{equation*}
This implies that the optimal 
solution on (\ref{op1r}) is unique, completing the proof of part (2) of the present {l}emma.
\endproof

\section{The optimal strategy, in continuum time}
\label{s:5}
\setcounter{equation}{0}
{W}e study 
the optimization problems {\bf (MTP)}, {\bf (OP)} on the entire plane,
assuming that the initial set 
$\Omega_0\subset \R^2$ is convex. 

According to Definition~\ref{d:11}, we say that the multifunction
$t\mapsto \Omega(t) \subset\R^2$  with compact values is {\bf admissible} if the following holds:

\begi
\item {\it 
For every time $t\geq 0$ and $\delta>0$,
\begin{equation}
\label{Equa:contiv_evol}
\Omega(t + \delta) \,\subset\, \Omega^{\delta} (t), \qquad \quad
\lim_{\delta  \to 0+} \frac{{\mathcal L^2(\Omega^\delta(t))} - {\mathcal L^2(\Omega(t+\delta))}}{\delta } \,= \, M.
\end{equation}
}
\endi
Assuming that all sets $\Omega(t)$ have perimeter with finite length, one has
\begin{equation*}
\begin{split}
\frac{d^+}{dt} {\mathcal L^2(\Omega(t))} ~&=~\lim_{\delta  \to 0+} \frac{{\mathcal L^2(\Omega(t+\delta))} - {\mathcal L^2(\Omega(t))}}{\delta } \\[4mm]
&\ds =~ \lim_{\delta  \to 0+} \frac{{\mathcal L^2(\Omega^\delta(t))} - {\mathcal L^2(\Omega(t))}}{\delta} - M ~
~= ~{\mathcal H^1(\partial \Omega(t))} - M.
\end{split}
\end{equation*}

For a given time interval $[0,T]$,
we are interested in finding the evolutions which minimize the terminal area
${\mathcal L^2(\Omega(T))}$.
Assuming that the initial set $\Omega(0)=\Omega_0$ is convex, we will show that these evolutions
also minimize the areas ${\mathcal L^2(\Omega(t))}$ of all intermediate sets, for $ t \in [0, T]$, respectively.

Let a bounded, open convex set $\Omega_0$ and a constant $M>0$  be given.
Using the notation introduced in  Lemma~\ref{l:31} and in Lemma~\ref{l:39}, we define an admissible strategy
by setting
\bel{Adef}
 A (t) ~=~ \Tilde\Omega \bigl(\Omega_0^t,a(t)\bigr),
\eeq
where the area function $a(t)= {\mathcal L^2(A(t))}$ satisfies
\bel{at} 
\frac{d}{dt} a(t) ~=~ {\mathcal H^1\big(\partial \Tilde\Omega \bigl(\Omega_0^t,a(t)\big)\big)} - M,\qquad\qquad a(0)\,=\, {\mathcal L^2(\Omega_0)}. 
\eeq
This strategy is defined on a time interval $t\in [0,T^*]$, where  
\bel{T*}
T^*~\doteq~\sup\bigl\{t>0\,;~a(t)>0\bigr\}~\in~\,]0, +\infty]\eeq
is the first time when the area vanishes. We will show that this is indeed the optimal strategy.

\begin{theorem}\label{t:61} Let a bounded, open convex set $\Omega_0$ and a constant $M>0$  be given.
Then the set valued map $A(\cdot)$ introduced at (\ref{Adef})--(\ref{T*}) 
is a well defined, admissible strategy.
For any other admissible strategy $\Omega(\cdot)$, at every time $t\in [0, T^*]$ the areas satisfy
\bel{comp}
 {\mathcal L^2(A(t))}~\leq~{\mathcal L^2(\Omega(t))}.
\eeq
As a consequence, one has
\begi
\item[(i)] Setting $A(t)=\emptyset$ for $t>T^*$, the map $t\mapsto A(t)$ provides a solution to the optimization problem {\bf (OP)}.  
\item[(ii)] The minimum time problem {\bf (MTP)} is solvable if and only if $T^*<+\infty$.
In this case, the map $A(\cdot)$ provides an optimal solution.
\endi
\end{theorem}

{\bf Proof.} {\bf 1.} We begin by showing that the function  $a(t)$ is well defined.
The proof of Theorem \ref{t:41} shows that the map
\begin{equation*}
(t,a)~ \mapsto~ {\mathcal H^1(\partial \Tilde\Omega (\Omega_0^t,a))}
\end{equation*}
is continuous w.r.t.~both variables and monotone increasing  w.r.t.~$a$.
It is also locally Lipschitz continuous w.r.t.~$a$ as long as
 $0 < {\mathcal L^2(\Tilde\Omega (\Omega_0^t,a))} < {\mathcal L^2(\Omega_0^t)}$.   
%
%
By Peano's theorem, a solution to  the ODE (\ref{at}) thus exists.

To prove uniqueness, consider 
 two solutions, say  $a(t) \leq a'(t)$. Observing that the maximum curvature of the 
 boundary of $\Omega_0^t$ is $\leq 1/t$, using (\ref{bcurv}) we obtain the bound
\begin{equation}
\label{Equa:leading_to_unique}
\begin{split}
\frac{d}{dt} \bigl(a'(t) - a(t)\bigr) ~&= {\mathcal H^1(\partial \Tilde \Omega(\Omega_0^t,a'(t))) - \mathcal H^1(\partial \Tilde \Omega(\Omega_0^t,a(t))) } \\[4mm]
\ds &{=} ~\int_{a(t)}^{a'(t)}
\frac{d}{da} {\mathcal H^1(\partial \Tilde\Omega (\Omega_0^t,a))} \, {da} ~ \leq~
\frac{a'(t) - a(t)}{t}\,,
\end{split}
\end{equation}
valid for $t>0$ sufficiently small.
Therefore {for $0 < s \leq t$}
\bel{uniq}
a'(t) - a(t) ~\leq ~\frac{t}{s} \bigl(a'(s) - a(s)\bigr).
\eeq
Since at $t=0$ one has
\begin{equation*}
\lim_{\delta \to 0+} \frac{{\mathcal L^2(\Omega_0^\delta)} - a(\delta)}{\delta}~ =~ \lim_{\delta \to 0+} \frac{{\mathcal L^2(\Omega_0^\delta)} - a'(\delta)}{\delta} ~=~ M,
\end{equation*}
we conclude that
\begin{equation*}
\lim_{s \to 0} \frac{a'(s) - a(s)}{s} ~=~ 0.
\end{equation*}
Letting $s\to 0$ in (\ref{uniq}),
this yields  the uniqueness of the solution. 

 By the 
definition of the area function at (\ref{at}), it follows that $t\mapsto A(t)$ is an admissible strategy.
\v
{\bf 3.} We now prove the inequalities (\ref{comp}), showing that the strategy (\ref{Adef}) is optimal.
Given any other admissible strategy $t\mapsto\Omega(t)$, consider the sets
\begin{equation*}
Z(t) ~\doteq~ \Tilde\Omega \Big(\Omega_0^t,\,{\mathcal L^2(\Omega(t))}\Big).
\end{equation*}
The time derivative of the function $z(t)\,\doteq\,{\mathcal L^2(Z(t))}\,=\,{\mathcal L^2(\Omega(t))}$ 
is computed by
$$\bega{rl}
\ds\frac{d^+}{dt}z(t) &\ds=~ \liminf_{\delta \to 0+} \frac{{\mathcal L^2(Z(t+\delta))} - {\mathcal L^2(Z(t))}}{\delta} \\[2mm]
&\ds =~ \lim_{\delta \to 0+} \frac{{\mathcal L^2(\Omega(t+\delta))} - {\mathcal L^2(\Omega^\delta(t))}}{\delta} + \liminf_{\delta \to 0+} \frac{{\mathcal L^2(\Omega^\delta(t))} - {\mathcal L^2(Z(t))}}{\delta} \\[2mm]
&\geq \,\ds- M + \liminf_{\delta \to 0+} \frac{|Z^{\delta}(t)| - {\mathcal L^2(Z(t))}}{\delta}\\[2mm]
&=~ |\partial Z(t)| - M.
\enda
$$
Therefore, $z(t) = {\mathcal L^2(\Omega(t))}$  satisfies the differential inequality
\begin{equation*}
\frac{d^+}{dt} z(t) ~\geq ~\bigl|\partial \Tilde\Omega (\Omega_0^t,z(t))\bigr| - M.
\end{equation*}
{Following the same line as in \eqref{Equa:leading_to_unique}, a} comparison with the optimal solution $a(t) = {\mathcal L^2(A(t))}$ now yields
$$
\frac{d^+}{dt} {(z(t) - a(t))} ~=~ \liminf_{\delta \to 0+} \frac{{(z(r+\delta)- a(t + \delta)) - (z(t) - a(t))}}{\delta} ~\geq~ \frac{{z(t)-a(t)}}{t},
$$ 
Since the
map  $t \mapsto {z(t)-a(t)}$ is continuous, we deduce 
\bel{zas}
{z(t)-a(t)~ \geq~ \frac{t}{s} (z(s)-a(s))}.
\eeq
On the other hand, since both strategies are admissible, at the initial time we have
\bel{laz}\lim_{\delta \to 0+} \frac{{z(\delta)-a(\delta)}}{\delta} ~=~ 0.
\eeq
Letting $s\to 0$ in (\ref{zas}) and using (\ref{laz}), for every $t\in [0, T^*]$ we thus obtain
\begin{equation*}
z(t)~ \geq ~a(t),
\end{equation*}
proving (\ref{comp}).

The two statements (i)-(ii)
are now an immediate consequence of (\ref{comp}).
\endproof

\begin{remark} {\rm The optimal strategy $A(t)$ introduced at (\ref{Adef})-(\ref{at}) 
is uniquely determined for all $t\in [0, T^*]$, thanks to the  formulas (\ref{TOM}) which remove any ambiguity also in the case (\ref{asmall}) where $a>0$ is very small.
In general, however, other optimal solutions can exist.

Calling $R(t)$ the inner radius of the set $\Omega_0^t$, 
if  $a(t)<\bigl| \Hat \Omega(\Omega^t_0, R(t))\bigr|$, according to
(iv) and (v) in Theorem~\ref{t:41}, 
the minimization problem
\bel{minp} \hbox{minimize:}\qquad {\mathcal H^1(\partial\Omega)}\qquad \hbox{subject to}\quad \Omega\subseteq
\Omega_0^t\,,\quad {\mathcal L}^2 (\Omega)= a(t),\eeq
has multiple solutions.  One can thus construct a different optimal strategy, say $t\mapsto A'(t)\subset\Omega_0^t$, where each set $A'(t)$ is a translation of the corresponding set $A(t)$.}
\end{remark}

\section{Large time behavior}
\label{s:6}
\setcounter{equation}{0}

In  this last section we study the large time behavior of the optimal strategy $A(t)$.   

\begin{proposition}\label{p:41} Let $\Omega_0\subset\R^2$ be a bounded, open  convex set.
Then there exists a constant $M_0>0$ such that the following holds.
\begi
\item[(i)] For $0<M<M_0$, the optimal strategy $A(t)$ defined at (\ref{Adef})--(\ref{T*}) satisfies
\bel{limi}\lim_{t\to +\infty} {\mathcal L^2(A(t))} ~=~ +\infty.\eeq
\item[(ii)] 
For $M=M_0$, after some time $T^\dagger\geq 0$ the set $A(t)$ becomes a ball:
\bel{Aco}
A(t)~=~B_{M/2\pi} (\bar x)\qquad \qquad\forall t\geq T^\dagger.\eeq
\item[(iii)]
 For $M>M_0$,  the area ${\mathcal L^2(A(t))}$ shrinks to zero  at a finite time $T^*<+\infty$.
In this case, there exists a time $T^\dagger<T^*$ such that the set $A(t)$ is a ball for all
$T^\dagger\leq t< T^*$.
\endi
\end{proposition}

{\bf Proof.}  {\bf 1.} Since our  solutions now depend on the choice of $M$, 
the notation $a_M(t) =
{\mathcal L^2(A_M(t))}$ will be used.
We consider the set of all solutions $a_M$ of (\ref{at}) which remain uniformly bounded
for all times $t\geq 0$, and define
\bel{M0def}
M_0~=~\inf\left\{M>0\,;~~\sup_{t>0}~ a_M(t) <+\infty\right\}.\eeq
{}From the equation (\ref{at}) it immediately follows
\bel{amm}
M'<M\qquad\implies\qquad  a_{M}(t)~\leq~a_{M'}(t)\qquad\qquad\forall t\geq 0.\eeq
\v
{\bf 2.} To prove (i), assume $M<M_0$.  By definition, $\sup_{t>0} a_M(t) = +\infty$.
By the isoperimetric inequality, any set $\Omega$ with area ${\mathcal L}^2 (\Omega)=a$ has perimeter 
${\mathcal H^1(\Omega)}\geq 2\sqrt{\pi a}$.
Hence from (\ref{at}) it follows
$${d\over dt} a_M(t)~\geq~ 2\sqrt{\pi a_M(t)} - M.$$
In particular,
if at some time $\tau>0$ one has  
$2\sqrt{\pi a_M(\tau)}  > M$, then the area function $t\mapsto a_M(t)$ is monotone increasing,
and approaches infinity as $t\to +\infty$. This proves part  (i).
\v
{\bf 3.} To prove (iii), assume $M>M_0$.   Let $M_0<M'<M$.   By assumption, the solution 
$a_{M'}(t)\geq a_M(t)$ is also uniformly bounded.  The difference between these two 
solutions of (\ref{at}) satisfies
$${d\over dt} \bigl( a_{M'}(t)-a_M(t)\bigr)~\geq~M-M'.$$
At all times $t\geq 0$ where $a_M(t)$ is defined, this implies
\bel{MM'} a_M(t)~\leq~\left(\sup_{\tau>0}~ a_{M'}(\tau) \right) - (M-M') t.\eeq
Since the right hand side of (\ref{MM'}) becomes negative for $t$ large, by continuity 
there exists a time $T^*$ such that ${\mathcal L^2(A_M(T^*))}= a_M(T^*)=0$.
This proves the first statement in (iii).  

Next,
call $R_0\geq 0$ the inner radius of the convex set $\Omega_0$, i.e.~the  radius of the 
largest ball contained inside $\Omega_0$. Then the inner radius of $\Omega_0^t$ is
$R(t)=R_0+t$.  By the definition (\ref{Adef}), when the area
satisfies $a_M(t)= {\mathcal L^2(A_M(t))}\leq \pi  R^2(t)$, the set $A_M(t)$ becomes a ball.
This is certainly true when $t$ is sufficiently close to $T^*$, because
$$\lim_{t\to T^*-} a_M(t)~=~0,\qquad\qquad \lim_{t\to T^*-} R(t)~=~R_0+T^*.$$
This establishes the last statement in (iii).
\v
{\bf 4.} Finally, consider the case $M=M_0$.    We claim that the corresponding solution
$t\mapsto a_M(t)$ remains uniformly positive, and uniformly bounded.

Indeed, assume that at some time $\tau>0$ one has 
$2\sqrt{\pi a_{M_0}(\tau)}  > M_0$.   By continuity, there exists $\ve>0$ such that 
$2\sqrt{\pi a_{M'}(\tau)}  > M'$ for every $M'\in [M_0, M_0+\ve]$ as well.  
As remarked in step {\bf 2}, this implies
$a_{M'}(t)  \to +\infty$ as $t\to +\infty$, reaching a contradiction with 
the definition of $M_0$.

On the other hand, assume that ${\displaystyle \liminf_{t\to\infty} a_{M_0}(t) = 0}$. In this case, following the argument in step {\bf 3}, we can find a time $\tau$ such that
$A_{M_0}(\tau)$ is a ball of radius $R<{M\over 2\pi}$.  By continuity there exists $\ve>0$ 
such that, for every $M'\in [M_0-\ve, M_0]$,  the set $A_{M'}(\tau)$ is also a ball of radius $R<{M'\over 2\pi}$.   If this happens, then the area $A_{M'}(\tau)$ shrinks to zero in finite time.
 Again, this yields  a contradiction with 
the definition of $M_0$.

It remains to prove that, for $t$ sufficiently large, $A_{M_0}(t)$ is a ball. Toward this goal we observe that the inner radius of $\Omega_0^t$ is 
$R(t)> t$.   Hence, when $a=a_{M_0}(t)< \pi t^2$, the solution to the optimization
problem (\ref{minp}) is a ball.   Since $a_{M_0}(t)$ remains bounded, this inequality is true (and hence $A_{M_0}(t)$ is a ball)
for all times $t$ large enough.   It is now clear that the radius of this ball must be
$R= {M\over 2\pi}$.  Otherwise, the solution of (\ref{at}) will tend to infinity, or else become zero 
in finite time. 
\endproof

\begin{remark} {\rm When the initial set $\Omega_0$ is a ball of radius $R$, one has
$M_0=2\pi R = 2 \sqrt{\pi {\mathcal L^2(\Omega_0)}}$.    
 For a general  convex set $\Omega_0$, one has 
$$M_0~\leq~2 \sqrt{\pi {\mathcal L^2(\Omega_0)}}\,.$$
Indeed, this threshold is smaller than the threshold for the ball, since the perimeter remains larger for the same area.

In the case where $\Omega_0$ is a square, for various values of $M$ the 
optimal strategy $A_M(t)$ has been 
studied  in Example 8.1 of \cite{BCS2}.
In particular, for the unit square the value of $M_0$ can be computed explicitly:
$$ M_0\,= \,{4-\pi \over 1-\ln 2}\,\approx \,2.797\,<\, 2\sqrt \pi.$$
In this case,  for $M=M_0$, the optimal set $A(t)$ becomes a ball at time 
$$T^\dagger={1 \over 4(1-\ln2)}-{1\over 2}.$$

}
\end{remark}

{\bf Acknowledgments.} The third author acknowledges the support from the NSERC Discovery Grant RGPIN-2023-04561.

\end{document}